\documentclass[twoside,12pt]{report}
\usepackage{graphicx}
\usepackage{amssymb, amsmath, amsxtra}
\usepackage{epsfig,subfigure}

\usepackage{latexsym}
\usepackage{bm}

 \usepackage{ntheorem}
           \theoremstyle{plain}
                      {\theorembodyfont{\rmfamily}
                      \theoremseparator{.}
                       \newtheorem{exa}{Example}[section]
           \newtheorem{thm}{Theorem}[section]
           \theoremstyle{plain}
           
           \theoremstyle{plain} \newtheorem{lem}{Lemma}[section]
           \theoremstyle{plain}

           \theoremstyle{plain}
           \newtheorem{rem}{Remark}[section]
            }

\begin{document}
\begin{center}
{\Large \textbf{MEMOIR ON INTEGRATION OF\\[.3ex] ORDINARY DIFFERENTIAL\\[1.2ex] EQUATIONS
BY QUADRATURE}}\footnote{First published in \textit{Archives of ALGA}, vol. 5, 2008, pp. 27-62.}\\[2ex]
 Nail H. Ibragimov\\
 Department of Mathematics and Science, Blekinge Institute
 of Technology,\\ 371 79 Karlskrona, Sweden
 \end{center}

 \noindent
\textbf{Abstract.} The Riccati equations
 reducible to first-order linear equations by an appropriate change
  the dependent variable are singled out. All these equations are
  integrable by quadrature.

  A wide class of linear ordinary differential equations
 reducible to algebraic equations is found. It depends on two arbitrary functions.
 The method for solving all these equations is given.
 The new class contains the constant coefficient equations and Euler's equations as
 particular cases.\\[2ex]

The old-fashioned title of the present paper indicates that it is
dedicated to quite old topics. Namely, it deals with a problem on
integration by quadrature of Riccati equations investigated, in
terms of elementary functions, by Francesco Riccati and Daniel
Bernoulli
 some 280 years ago for the special Riccati equations (see, e.g.
 \cite{step58})
 $$
 y' = a y^2 + b x^\alpha,\quad a, b, \alpha = {\rm const.},
 $$
  and with an integration of higher-order linear equations by reducing
  them to algebraic equations. The later property was discovered
  by Leonard Euler in the 1740s for the constant coefficient
  equations
 $$
 y^{(n)} + A_1 y^{(n-1)} + \cdots + A_{n-1} y' + A_n y = 0,
 \quad A_1, \ldots, A_n = {\rm const.,}
 $$
 as well as for the equations of the form
 $$
 x^ny^{(n)} +  A_1 x^{n-1} y^{(n-1)} + \cdots + A_{n-1} x y' + A_n y = 0,
\quad A_1, \ldots, A_n = {\rm const.,}
 $$
known as Euler's equations.

 It will be shown in what follows that these classical results can be extended to wide classes of equations.

 \chapter{Riccati equation}
% \addcontentsline{toc}{chapter}{Riccati equation}

 \section{Introduction}
% \addcontentsline{toc}{section}{Introduction}
 \label{Mint.int}
 \setcounter{equation}{0}

Consider the special Riccati equation
 \begin{equation}
 \label{Ricc.special}
 y' = a y^2 + b x^\alpha, \quad a, b, \alpha = {\rm const.}
 \end{equation}
 If $\alpha = 0,$ Eq. (\ref{Ricc.special}) is integrable by the separation of variables:
 $$
 \frac{dy}{a y^2 + b} = dx.
 $$
 Another easily integrable case is $\alpha = - 2.$ Then the change of the dependent variable
 $$
 z = \frac{1}{y}
 $$
 maps Eq. (\ref{Ricc.special}) to the homogeneous equation
 $$
 \frac{dz}{dx} = - \Big[a - b \Big(\frac{z}{x}\Big)^2\Big]
 $$
 which is integrable by quadrature.

 F. Riccati and D. Bernoulli noted independently that Eq. (\ref{Ricc.special})
 can be transformed to the case $\alpha = 0,$ and hence integrable by quadrature
 in terms of elementary functions if $\alpha$ takes the values from the following
 two series:
 \begin{align}
 \label{Ricc.special_1}
  & \alpha  = - 4, \ - \frac{8}{3}\,, \ - \frac{12}{5}\,,
 \ - \frac{16}{7}\,,  \ldots\,;\notag\\[1ex]
   & \alpha  = - \frac{4}{3}\,, \ - \frac{8}{5}\,,
 \ - \frac{12}{7}\,, \ - \frac{16}{9}\,, \ldots\,.
 \end{align}
The series (\ref{Ricc.special_1}) are given by the formula
 \begin{equation}
 \label{Ricc.elementary}
  \alpha = - \frac{4k}{2k - 1} \quad  {\rm with} \quad k =  \pm 1, \ \pm
 2, \ldots\, .
 \end{equation}
It is manifest from (\ref{Ricc.elementary}) that both series in
(\ref{Ricc.special_1}) have the limit $\alpha = - 2.$ For a
derivation of the transformations mapping Eq. (\ref{Ricc.special})
with $\alpha$ having the form (\ref{Ricc.elementary}) to an
integrable form, see \cite{step58}, Chapter 1, \S 6.

 J. Liouville showed in 1841
 that the solution to the special Riccati equation
 (\ref{Ricc.special}) is integrable by quadrature
 in terms of elementary functions only if $\alpha$ has the
 form (\ref{Ricc.elementary}).

 It is well known that the general Riccati equation
 $$
 y' = P(x) + Q(x) y + R(x) y^2
 $$
 can be rewritten as a linear second-order equation. But this kind of linearization
 by \textit{raising the order} does not solve the integration problem.
 Therefore, I will investigate a possibility of linearization of
 the Riccati equations \textit{without raising the order} and will show that all
 Riccati equations of this type can be integrated by quadrature.

 \section{The linearizable Riccati equations}
 \label{gai.sub1}
  \setcounter{equation}{0}

 The following theorem is closely related to the theory of
 nonlinear superpositions discussed in \cite{ibr06P}, Section 6.7.

  \begin{thm}
 \label{gai.th0}
The first-order ordinary differential equation
 \begin{equation}
 \label{gai.eq1}
 y' = f(x, y)
 \end{equation}
 can be reduced to a linear first-order equation
 \begin{equation}
 \label{gai.eq6}
 \frac{d z}{dx} = p(x) + q(x) z
 \end{equation}
  by a change  of the dependent variable $y,$
 \begin{equation}
 \label{gai.th.eq}
 z = z(y),
 \end{equation}
  if and only if
 Eq. (\ref{gai.eq1}) can be written in the form
 \begin{equation}
 \label{gai.th.eq4}
 y' = T_1(x) \xi_1(y) + T_2(x)  \xi_2(y)
 \end{equation}
 such that the operators
 \begin{equation}
 \label{gai.th.eq5}
 X_1 = \xi_1(y)\frac{\partial}{\partial y}\,, \quad
 X_2 = \xi_2(y) \frac{\partial}{\partial y}
 \end{equation}
 span a two-dimensional (or a one-dimensional if $X_1$ and
 $X_2$ are linearly dependent) Lie algebra called in \cite{ibr99} the
 VGL (Vessiot-Guldberg-Lie) algebra.
  \end{thm}
 {\bf Proof.}  Let Eq. (\ref{gai.eq1}) be linearizable.  Then we can assume that it is already reduced by a certain change of
 the dependent variable
  \begin{equation}
 \label{gai.th0.eq}
  z = \zeta(y)
 \end{equation}
  to a linear equation (\ref{gai.eq6}),
 $$
 \frac{d z}{dx} = p(x) + q(x) z.
 $$
 The VGL algebra associated with Eq.
(\ref{gai.eq6}) is two-dimensional and is spanned by the operators
 \begin{equation}
 \label{gai.eq6vgl}
 \overline X_1 = \frac{\partial}{\partial z}\,, \quad
 \overline X_2 = z\frac{\partial}{\partial z}\,\cdot
 \end{equation}
The form of Eq. (\ref{gai.th.eq4})
 and the algebra property $[X_1, X_2] = \alpha X_1 + \beta
 X_2$ remain unaltered under any change (\ref{gai.th.eq}) of the
 dependent variable. Therefore, rewriting Eq. (\ref{gai.eq6}) in
 the original variable $y = \zeta^{-1}(z)$ obtained by the inverse transformation
 to (\ref{gai.th0.eq}), we arrive at an equation of the form (\ref{gai.th.eq4})
 for which the operators (\ref{gai.th.eq5}) span a two-dimensional
 Lie algebra. Since the equation obtained from Eq. (\ref{gai.eq6})
  by the inverse transformation to (\ref{gai.th0.eq}) is the
  original equation (\ref{gai.eq1}) we have proved the ``only if" part
  of the theorem.

  Let us prove now the ``if" part  of the theorem. Namely, we have to demonstrate that any
  equation of the form (\ref{gai.th.eq4}) such that the operators
   (\ref{gai.th.eq5}) span a two-dimensional Lie algebra, is linearizable.
 If the operators
 (\ref{gai.th.eq5}) are linearly dependent, then  $\xi_2(x) =
 \gamma \xi_1(x), ~\gamma = {\rm const.,}$ and hence Eq.
(\ref{gai.th.eq4}) has the form
$$
 y' = \big[T_1(x) + \gamma T_2(x)\big] \xi_1(y).
$$
It can be reduced to the linear equation
 $$
 z' = T_1(x) + \gamma T_2(x).
 $$
 upon introducing a canonical variable $z$ for
 $$
 X_1 = \xi_1(y) \frac{\partial}{\partial y}
 $$
 by solving the equation $X_1(z)= 1.$

Suppose now that the operators (\ref{gai.th.eq5}) are linearly
independent. It is clear that Eq. (\ref{gai.ric2}) will be
linearized if one transforms the operators (\ref{gai.th.eq5}) to the
form  (\ref{gai.eq6vgl}). One can assume that the first operator
(\ref{gai.th.eq5}) has been already written, in a proper variable
$z,$ in the form of the first operator $\overline X_1$ given in
(\ref{gai.eq6vgl}):
 $$
 X_1 = \frac{\partial}{\partial z}\,\cdot
 $$
 Let the second operator
(\ref{gai.th.eq5}) be written in the variable $z$ as follows:
 $$
 X_2 = f(z)\frac{\partial}{\partial z}\,\cdot
 $$
We have
 $$
 [X_1, X_2] = f'(z) \frac{\partial}{\partial z}
 $$
 and the requirement  $[X_1, X_2] = \alpha X_1
 + \beta X_2$ that $X_1, X_2$ span a Lie algebra  $L_2$ yields the differential equation
 $$
 f' = \alpha + \beta f,
 $$
where not both  $\alpha$ and $\beta$ vanish because otherwise  the
operators $X_1$ and $X_2$ will be linearly dependent. hence $f'(x)
\not= 0.$ Solving the above differential equation, we obtain
 \begin{align}
 & f = \alpha z + C \quad {\rm if} \quad \beta = 0,\notag\\
 & f = C\,{\rm e}^{\beta z} - \frac{\alpha}{\beta} \quad
{\rm if} \quad \beta \not= 0.\notag
 \end{align}
$$
f = a x + C \quad \Longrightarrow \quad X_2 = ax \frac{{\rm d}}{{\rm
d} x} + C X_1, \quad {\rm if} \quad b = 0,
$$
$$
f = C{\rm e}^{bx} - \frac{a}{b} \quad \Longrightarrow \quad X_2 =
 C{\rm e}^{bx} \frac{{\rm d}}{{\rm d} x} - \frac{a}{b} X_1,\quad
{\rm if} \quad b \not= 0.
$$
In the first case we have
 $$
 X_2 = \alpha z \frac{\partial}{\partial z} + C X_1,
 $$
 and hence a basis of $L_2$ is provided by (\ref{gai.eq6vgl}).
 In the second case we have
 $$
 X_2 = C\,{\rm e}^{\beta z}\, \frac{\partial}{\partial z}
 - \frac{\alpha}{\beta}\, X_1,
 $$
and one can take basis operators, by assigning $\beta z$ as new $z,$
in the form
 $$
 X_1 = \frac{\partial}{\partial z}\,, \quad
  X_2 = {\rm e}^z\, \frac{\partial}{\partial z}\,\cdot
 $$
 Finally, substituting $\overline z = {\rm e}^{-z}$ we arrive at
 the basis (\ref{gai.eq6vgl}), thus completing the proof of the theorem.

 The following theorem characterizes all Riccati equations that can be
 reduced to first-order linear equations by changing
  the dependent variable (see \cite{ibr89}, Russian ed., Theorem 4.3; \cite{ibr99}, Section 11.2.5 and
  Note [11.4]; see also Theorem 3.2.2 in
  \cite{ibr06P}).

  \begin{thm}
 \label{gai.th}
The following two conditions $1^\circ$ and $2^\circ$ are equivalent
and provide the necessary and sufficient conditions for the Riccati
equation
 \begin{equation}
 \label{gai.ric1}
 y' = P(x) + Q(x) y + R(x) y^2
 \end{equation}
 to be linearizable by a change of the dependent variable (\ref{gai.th.eq}),
 $z = z(y).$

 $\bm{1^\circ.}$
 Eq. (\ref{gai.ric1}) has a constant solution $y = c$ including $c =
 \infty.$

 $\bm{2^\circ.}$  Eq. (\ref{gai.ric1}) has either the form
 \begin{equation}
 \label{gai.ric2}
 y' = Q(x)y + R(x)y^2
 \end{equation}
 with any functions $Q(x)$ and $R(x),$ or the form
 \begin{equation}
 \label{gai.ric3}
 y' = P(x) + Q(x)y + k[Q(x)- kP(x)]y^2
 \end{equation}
 with any functions $P(x), ~Q(x)$ and any constant $k.$
 \end{thm}
 {\bf Proof.} The conditions $1^\circ$ and $2^\circ$ are equivalent.
 Indeed, let Eq. (\ref{gai.ric1}) have a constant solution $y = c.$
 Then
 \begin{equation}
 \label{gai.ric4}
 0 = P(x) + Q(x) c + R(x) c^2.
 \end{equation}
 If $c = 0,$ Eq. (\ref{gai.ric4}) yields $P(x) = 0,$ and
 hence Eq. (\ref{gai.ric1}) has the form (\ref{gai.ric2}).
 If $c \not= 0,$ Eq. (\ref{gai.ric4}) yields
 $$
 R(x) = - \frac{1}{c^2} \big[P(x) + c Q(x)\big] = - \frac{1}{c} \big[Q(x) + \frac{1}{c}\,
 P(x)\big].
 $$
Denoting $k = - 1/c$ we get
 $$
 R(x) = k \big[Q(x) -k P(x)\big].
 $$
 Hence Eq. (\ref{gai.ric1}) has the form (\ref{gai.ric3}). Thus,
 we have proved that $1^\circ \Rightarrow 2^\circ.$

 Conversely, let Eq. (\ref{gai.ric1}) satisfy the condition $2^\circ.$
 It is manifest that Eq. (\ref{gai.ric2}) has the constant solution $y =
 0.$ Furthermore, one can verify that  Eq. (\ref{gai.ric3}) has
 the constant solution $y = - 1/k.$ This proves that $2^\circ \Rightarrow 1^\circ.$

 Let us turn to the necessary and sufficient conditions for
 linearization. The operators
 $$
  X_1 = y \frac{\partial}{\partial y}\,, \quad
  X_2 = y^2\frac{\partial}{\partial y}
 $$
 associated with Eq. (\ref{gai.ric2}) have the commutator $[X_1, X_2] = X_2,$
 and hence span a two-dimensional Lie algebra. Furthermore, the operators
 $$
  X_1 = (1 - k^2y^2) \frac{\partial}{\partial y}\,, \quad X_2 =
 (y + k y^2)\frac{\partial}{\partial y}
 $$
 associated with Eq. (\ref{gai.ric3}) have the commutator $[X_1, X_2] = X_1 + 2 k X_2.$
 Hence, they also span a two-dimensional Lie algebra.  Therefore,
 according to Theorem \ref{gai.th0}, the condition $2^\circ$
 is sufficient for linearization.

 Note that the  equation $y' = P(x) + Q(x)y$ can be regarded as a
 particular case of Eq. (\ref{gai.ric3}) with $k = 0.$ Since Eq. (\ref{gai.ric3})
 has the constant solution $y = - 1/k,$ we conclude that the
 linear equation $y' = P(x) + Q(x)y$ has the constant
 solution $y = \infty.$ We conclude that any
 linearizable equation has a constant solution because the change of the dependent
 variable $z = z(y)$ maps a constant solution into a constant
 solution of the transformed equation. Hence, the condition $1^\circ$
 is necessary  for linearization. This completes the proof of the theorem
 due to the equivalence of the conditions $1^\circ$ and $2^\circ.$

 Now we will  use Theorem \ref{gai.th}
 for integrating the linearizable Riccati equations.
 We will find linearizing transformations
 for the equations (\ref{gai.ric2}) and (\ref{gai.ric3}).
 We will assume that $k$ in Eq. (\ref{gai.ric3}) is a real number.

 \section{Linearization and integration of  Eq.
 (\ref{gai.ric2})}
 \label{gai.sub2}
 \setcounter{equation}{0}

 Invoking Eqs. (6.7.5), (6.7.6) from \cite{ibr06P}, replacing there
  $t$ and $x^i$ by $x$ and $y,$ respectively, and identifying $T_1(t)$ and $T_2(t)$
 with $Q(x)$ and $R(x),$ respectively,  we see
 that the VGL (Vessiot-Guldberg-Lie) algebra
% Vessiot-Guldberg-Lie algebra (VGL algebra)
 associated with Eq. (\ref{gai.ric2}) is a two-dimensional
 algebra $L_2$ spanned by the operators
 $$
  X_1 = y \frac{\partial}{\partial y}\,, \quad X_2 =
  y^2\frac{\partial}{\partial y}\,\cdot
 $$
 Their commutator is $[X_1, X_2] = X_2.$ Introducing
  the new basis
  $$
  X'_1 = X_2, \quad X'_2 = - X_1
  $$
  i.e. taking
 \begin{equation}
 \label{gai.eq4}
  X'_1 = y^2 \frac{\partial}{\partial y}\,, \quad X'_2 =
  - y\frac{\partial}{\partial y}
 \end{equation}
 we have a basis in $L_2$
 satisfying the commutator relation
 \begin{equation}
 \label{gai.eq5}
 [X'_1, X'_2] = X'_1.
 \end{equation}
Let us  find a change of the dependent variable, $z = z(y),$ such
 that  Eq. (\ref{gai.ric2}) becomes a linear equation
 (\ref{gai.eq6}),
 $$
 \frac{d z}{dx} = p(x) + q(x) z.
 $$
The VGL algebra of this equation is spanned by the operators
(\ref{gai.eq6vgl}),
 $$
 \overline X_1 = \frac{\partial}{\partial z}\,, \quad
 \overline X_2 = z\frac{\partial}{\partial z}
 $$
 whose commutator has the same form as Eq. (\ref{gai.eq5}), i.e.
 $[\overline X_1, \overline X_2] = \overline X_1.$ Consequently,
 the linearizing transformation  $z = z(y)$ is determined by
 the equations $X'_1(z) = 1, \ X'_2(z) = z,$ or
 \begin{equation}
 \label{gai.eq7}
  y^2 \frac{d z}{dy} = 1, \quad - y \frac{d z}{dy} = z.
 \end{equation}
 Integrating the first equation (\ref{gai.eq7}) we obtain
 $$
 z = - \frac{1}{y} + A.
 $$
 Substituting in the second equation (\ref{gai.eq7}) we get
 $A = 0.$ Thus,  the linearizing transformation is
 \begin{equation}
 \label{gai.eq8}
 z = - \frac{1}{y}\,\cdot
 \end{equation}
 In this variable, Eq. (\ref{gai.ric2}) becomes the
 linear equation
 \begin{equation}
 \label{gai.eq9}
 z' = R(x) - Q(x) z,
 \end{equation}
whence
 $$
 z = \Big[C + \int R(x) {\rm e}^{\int Q(x) dx}\,dx\Big]\,
 {\rm e}^{-\int Q(x) dx}, \quad C = {\rm const.}
 $$
Substituting in Eq. (\ref{gai.eq8}) we finally arrive at the
 solution to Eq. (\ref{gai.ric2}):
 \begin{equation}
 \label{gai.eq10}
 y = - \Big[C + \int R(x) {\rm e}^{\int Q(x) dx}\,dx\Big]^{-1}\,
 {\rm e}^{\int Q(x) dx}.
 \end{equation}
 \begin{rem}
 In Section \ref{Ricc.Lin_subsec1} the solution (\ref{gai.eq10}) is obtained by an alternative method.
 \end{rem}

 \begin{exa}
 Consider the equation
 $$
 y'  = \frac{y}{x} + R(x) y^2.
 $$
 By Eq. (\ref{gai.eq10}) yields its solution
$$
 y = \frac{x}{C - \int x R(x) dx}\,\cdot
 $$
 \end{exa}

 \begin{exa}
 The equation
 $$
 y'  = \frac{y}{x} + 3 x y^2
 $$
is a particular case of the previous equation. Its solution is
$$
 y = \frac{x}{C - x^2}\,\cdot
 $$
 \end{exa}

 \begin{exa}
 The equation
 $$
 y'  = y + \frac{1}{x}\, y^2
 $$
has the solution
$$
 y = \frac{{\rm e}^x}{C + {\rm Ei}(x)}\,, \quad {\rm where}
 \quad {\rm Ei}(x) = -  \int\limits_{- \infty}^x \frac{{\rm e}^t}{t}\,dt.
 $$
 \end{exa}

 \section{Linearization and integration of  Eq.
 (\ref{gai.ric3})}
 \label{gai.sub3}
 \setcounter{equation}{0}

Now we identify the coefficients $P(x)$ and $Q(x)$ of Eq.
(\ref{gai.ric3}),
 \begin{equation}
 \tag*{(\ref{gai.ric3})}
 y' = P(x) + Q(x)y + k[Q(x)- kP(x)]y^2,
 \end{equation}
 with the coefficients $T_1(t)$ and $T_2(t)$
 of Eqs. (6.7.5), (6.7.6) from \cite{ibr06P} and associated with
 Eq. (\ref{gai.ric3}) the VGL algebra spanned by the operators
 $$
  X_1 = (1 - k^2y^2) \frac{\partial}{\partial y}\,, \quad X_2 =
 (y + k y^2)\frac{\partial}{\partial y}\,\cdot
 $$
 Their commutator is $[X_1, X_2] = X_1 + 2 k X_2.$ We assume that
 $k \not= 0,$ since otherwise Eq. (\ref{gai.ric3}) is already
 linear. Therefore, setting
 $$
 X'_1 = X_1 + 2 k X_2, \quad X'_2 = - X_1/(2k)
 $$
 we obtain the new basis
 \begin{equation}
 \label{gai.eq11}
  X'_1 = (1 + k y)^2 \frac{\partial}{\partial y}\,,
  \quad X'_2 = \frac{1}{2k}\, \big(k^2 y^2 - 1\big) \frac{\partial}{\partial y}\,,
 \end{equation}
satisfying the commutator relation (\ref{gai.eq5}),
 $[X'_1, X'_2] = X'_1.$

 The transformation $z = z(y)$ of the operators (\ref{gai.eq11})
 to the form (\ref{gai.eq6vgl})  is determined by
 the equations
 $$X'_1(z) = 1, \quad X'_2(z) = z,
 $$
 or
 \begin{equation}
 \label{gai.eq12}
  (1 + k y)^2 \frac{d z}{dy} = 1,
  \quad \frac{1}{2k}\, \big(k^2 y^2 - 1\big) \frac{d z}{dy} = z.
 \end{equation}
 Integrating the first equation (\ref{gai.eq12})
  we obtain
 $$
 z = A - \frac{1}{k (1 + k y)}\,\cdot
 $$
 Substituting in the second equation (\ref{gai.eq12}) we get
 $A = 1/(2k).$ Thus,  the linearizing transformation is
 \begin{equation}
 \label{gai.eq13}
 z = \frac{ky - 1}{2 k (k y + 1)}\,\cdot
 \end{equation}
 Eq. (\ref{gai.eq13}) yields:
 \begin{equation}
 \label{gai.eq14}
 y = \frac{1 + 2 k z}{k (1 - 2 k z)}\,\cdot
 \end{equation}
 Substituting (\ref{gai.eq14}) in Eq. (\ref{gai.ric3}) we arrive at
 the linear equation
 \begin{equation}
 \label{gai.eq15}
 z' = \frac{1}{2 k} \, Q(x) + \big[Q(x) - 2 k P(x)\big] z,
 \end{equation}
 whence
 \begin{equation}
 \label{gai.eq16}
 z = \frac{1}{2k} \left[C + \int Q(x) {\rm e}^{\int [2 k P(x) - Q(x)] dx}\, dx\right]\,
 {\rm e}^{\int [Q(x) - 2k P(x)] dx}.
 \end{equation}
Substituting (\ref{gai.eq16}) in  (\ref{gai.eq14}) we will obtain
the solution to Eq. (\ref{gai.ric3}).
 \begin{rem}
 \label{Ric3:rem1}
 We have assumed that $Q(x) \not= 0$ because otherwise Eq. (\ref{gai.ric3}) is separable,
 $$
 y' = P(x)(1 - k^2y^2).
 $$
 \end{rem}
 \begin{rem}
 \label{Ric3:rem2}
 Eq. (\ref{gai.ric3}) has the constant solution
 \begin{equation}
 \label{Ric3:sing}
 y = - \frac{1}{k}\,\cdot
 \end{equation}
 We excluded the singular solution (\ref{Ric3:sing}) in the
 above calculations and assumed that $1 + k y \not= 0.$
 \end{rem}

 \begin{exa}
 \label{gai.sub2.ex1}
 Let us integrate the equation
 \begin{equation}
 \label{gai.eq17}
 y' = x + 2 xy + x y^2.
 \end{equation}
 Here
 $$
 P(x) = x,\quad Q(x) = 2 x,\quad k = 1.
 $$
 Therefore the linearized equation (\ref{gai.eq15}) is written
 $z' = x.$ Integrating it and denoting the constant of integration by $C/2,$
 we obtain
 $$
 z =\frac{1}{2}\,(C + x^2)
 $$
 in accordance with  Eq. (\ref{gai.eq16}). Substituting the above $z$ in
 Eq. (\ref{gai.eq14}) we obtain the following solution to
 Eq. (\ref{gai.eq17}):
 \begin{equation}
 \label{gai.eq17sol}
 y = \frac{1 + C + x^2}{1 - C - x^2}\,\cdot
 \end{equation}
 Remark \ref{Ric3:rem2} gives also the singular solution $y = - 1.$
 \end{exa}

 \begin{rem}
 Eq. (\ref{gai.eq17}) can also be integrated by separating the
 variables.
 \end{rem}

  \begin{exa}
 \label{gai.sub2.ex2}
 The equation
 \begin{equation}
 \label{gai.eq18}
 y' = x^2 + (x + x^2)\, y + \frac{1}{4} \big(2 x +
 x^2\big)\, y^2
 \end{equation}
 has the form (\ref{gai.ric3}) with
 $$
 P(x) = x^2, \quad Q(x) = x + x^2, \quad k = \frac{1}{2}\,\cdot
 $$
 The linearized equation (\ref{gai.eq15}) is
written
 $$
 z' = x + x^2  + x z
 $$
 and has the solution
 $$
 z = \left(C + \int (x + x^2) {\rm e}^{- x^2/2}\,dx \right) {\rm
 e}^{x^2/2}.
 $$
 Substitution  in Eq. (\ref{gai.eq14}) yields
 $$
 y = 2 \frac{1 + z}{1 - z}\,\cdot
 $$
 Hence, the solution to Eq. (\ref{gai.eq18}) is given by
 \begin{equation}
 \label{gai.eq19}
 y = 2 \frac{1 + \left(C + \int (x + x^2) {\rm e}^{- x^2/2}\,dx \right) {\rm
 e}^{x^2/2}}{1 - \left(C + \int (x + x^2) {\rm e}^{- x^2/2}\,dx \right) {\rm
 e}^{x^2/2}}\,\cdot
 \end{equation}
  \end{exa}

 \chapter{Higher-order linear equations reducible to algebraic equation}
   \setcounter{section}{4}
  \section{Introduction}
 \label{NC.int}
 \setcounter{equation}{0}

 The linear ordinary differential equations with constant
coefficients
 \begin{equation}
 \label{NC.int1}
 y^{(n)} + A_1 y^{(n-1)} + \cdots + A_{n-1} y' + A_n y = 0,
 \quad A_1, \ldots, A_n = {\rm const.,}
 \end{equation}
 and Euler's equations
  \begin{equation}
 \label{NC.int2}
 x^ny^{(n)} +  A_1 x^{n-1} y^{(n-1)} + \cdots + A_{n-1} x y' + A_n y = 0,
\quad A_1, \ldots, A_n = {\rm const.,}
 \end{equation}
 are discussed practically in all textbooks on differential equations.
 They are useful in applications. The most remarkable property is that
 both equations (\ref{NC.int1}) and (\ref{NC.int2}) are \textit{reducible
 to algebraic equations.} Namely,
 their fundamental systems of solutions, and hence the general solutions
 can be obtained by solving algebraic equations.
 Then one can integrate the corresponding non-homogeneous linear equations by
 using the method of variation of parameters.
 It is significant to understand the nature of the reducibility
 and to extend the class of linear
 equations reducible to algebraic equations.

 In the present paper we will investigate this problem and find a wide class of
 linear ordinary differential equations
 that are reducible to algebraic equations.  The new class depends on two arbitrary
 functions of $x$ and contains the equations (\ref{NC.int1}) and
 (\ref{NC.int2}) as particular cases. The method for solving all
 these equations is given in Sections \ref{NC.gen} and \ref{NC.gen_third},
  and illustrated by examples in Section \ref{NC.gen_exs}.

  The main statement for the second-order equations is as follows
  (Section \ref{NC.gen}).

 \begin{thm}
 \label{NC.int_theorem}
 The linear second-order equations
 $$
 P(x) y'' + Q(x) y' + R(x) y = F(x)
 $$
  whose general solution can be obtained by solving
 algebraic equations and by quadratures,
 have the form
 \begin{equation}
 \label{NC.int3nh}
 \phi^2\,y'' + \big(A + \phi' - 2 \sigma\big) \phi\, y' +
  \big(B - A \sigma + \sigma^2 - \phi \sigma'\big) y = F(x),
 \end{equation}
  where $\phi = \phi(x), \ \sigma = \sigma(x)$ and $F(x)$ are arbitrary
 (smooth) functions, and $A, B = {\rm const.}$ The homogeneous
equation (\ref{NC.int3nh}),
 \begin{equation}
 \label{NC.int3}
 \phi^2\,y'' + \big(A + \phi' - 2 \sigma\big) \phi\, y' +
  \big(B - A \sigma + \sigma^2 - \phi \sigma'\big) y = 0,
 \end{equation}
 has the solutions of the form
 \begin{equation}
 \label{NC.int4}
 y = {\rm e}^{\int \frac{\sigma (x) + \lambda}{\phi(x)}\,dx},
 \end{equation}
 where $\lambda$ satisfies the \textit{characteristic equation}
 \begin{equation}
 \label{NC.int5}
 \lambda^2 + A \lambda + B = 0.
 \end{equation}
If the characteristic equation (\ref{NC.int5}) has
 distinct real roots $\lambda_1 \not= \lambda_2,$ the general solution
 to  Eq. (\ref{NC.int3}) is given by
 \begin{equation}
 \label{NC.int6}
 y(x) = K_1\,{\rm e}^{\int \frac{\sigma (x) + \lambda_1}{\phi(x)}\,dx}
  + K_2\,{\rm e}^{\int \frac{\sigma (x) +
  \lambda_2}{\phi(x)}\,dx}, \quad K_1, K_2 = {\rm const.}
 \end{equation}
 In the case of complex roots, $\lambda_1 = \gamma + i \theta, \
 \lambda_2 = \gamma - i \theta,$
 the general solution to  Eq. (\ref{NC.int3}) is given by
 \begin{equation}
 \label{NC.int7}
 y(x) = \left[K_1\,\cos \left(\theta \int \frac{dx}{\phi (x)}\right)
 + K_2\,\sin \left(\theta \int \frac{dx}{\phi (x)}\right)\right]
  {\rm e}^{\int \frac{\sigma (x) + \gamma}{\phi(x)}\,dx}.
 \end{equation}
 If the characteristic equation (\ref{NC.int5}) has equal roots
 $\lambda_1 = \lambda_2,$  the general solution to  Eq. (\ref{NC.int3})
 is given by
 \begin{equation}
 \label{NC.int8}
 y = \left[K_1 + K_2 \, \int \frac{dx}{\phi(x)}\right]
 {\rm e}^{\int \frac{\sigma (x) + \lambda_1}{\phi(x)}\,dx}, \quad
 K_1, K_2 = {\rm const.}
 \end{equation}
 The general solution of the non-homogeneous
equation (\ref{NC.int3nh}) can be obtained by the method of
variation of parameters.
 \end{thm}

\begin{rem}
 The equations with constant coefficients and Euler's equation are
 the simplest representatives of Eqs. (\ref{NC.int3}). Namely,
 setting $\phi (x) = 1, \ \sigma (x) = 0$ we obtain the second-order equation with constant
 coefficients
 $$
 y'' + A y' + B y = 0,
 $$
 and Eq. (\ref{NC.int4}) yields the well-known formula
 $$
 y = {\rm e}^{\lambda \, x},
 $$
 where $\lambda$ is determined by the characteristic equation
 (\ref{NC.int4}).

 If we set $\phi (x) = x, \ \sigma (x) = 0,$ we obtain the second-order Euler's
 equation (\ref{NC.int2}) written in the form
 $$
 x^2 y'' + (A + 1) x y' + B y = 0.
 $$
 Then Eq. (\ref{NC.int4}) yields the particular solutions for
 Euler's equations:
 $$
 y = x^\lambda,
 $$
 where $\lambda$ is determined  again by the characteristic equation
 (\ref{NC.int4}). For details, see Section \ref{NC.gen_exs}.
 For other functions  $\phi (x)$ and $\sigma (x) = 0,$ Eqs.
 (\ref{NC.int4}) are new.
\end{rem}

 \section{Constant coefficient and Euler's equations from the group standpoint}
 \label{NC.gp}
 \setcounter{equation}{0}

For the sake of simplicity, we will consider in this section
second-order equations
 \begin{equation}
 \label{NC.eq1}
 y'' + f(x) y' + g(x) y = 0.
 \end{equation}

 \subsection{Equations with constant coefficients}
 \label{NC.cc}

 Let us begin with the equations with constant coefficients
 \begin{equation}
 \label{NC.eq2}
  y'' + A y' + B y = 0, \quad A, B = {\rm const.}
 \end{equation}
 Eq. (\ref{NC.eq2}) is invariant under
 the one-parameter groups of translations in $x$ and dilations in $y,$ since
 it does not involve the independent variable $x$ explicitly
 (the coefficients $A$ and $B$ are constant) and is homogeneous in the
 dependent variable $y.$
 In other words, Eq. (\ref{NC.eq2}) admits the generators
 \begin{equation}
 \label{NC.eq3}
 X_1 = \frac{\partial}{\partial x}\,, \quad  X_2 = y \frac{\partial}{\partial y}
 \end{equation}
of the translations in $x$ and dilations in $y.$ We use them as
follows \cite{ibr99}.

Let us find the invariant solution for $X = X_1 + \lambda X_2,$ i.e.
 \begin{equation}
 \label{NC.eq4}
 X =\frac{\partial}{\partial x} +
 \lambda y \frac{\partial}{\partial y}\, \quad \lambda = {\rm
 const.}
 \end{equation}
 The characteristic equation
 $$
 \frac{d y}{y} = \lambda \, \frac{d x}{x}
 $$
 of the equation
 $$
 X(J) \equiv \frac{\partial J}{\partial x} +
 \lambda y \frac{\partial J}{\partial y} = 0
 $$
 for the invariants $J(x, y)$ yields
 one functionally independent invariant
 $$
 J = y\,{\rm e}^{-\lambda x}.
 $$
 According to the general theory, the invariant solution is given by $J= C$
with an arbitrary constant $C.$ Thus, the general form of the
invariant solutions for the operator (\ref{NC.eq4}) is
 $$
 y = C {\rm e}^{\lambda x}, \quad C = {\rm const.}
 $$
 Since Eq. (\ref{NC.eq2}) is homogeneous one can set $C = 1$ and obtain \textit{Euler's
substitution}:
 \begin{equation}
 \label{NC.eq5}
  y = {\rm e}^{\lambda x}.
  \end{equation}
 As well known, the substitution (\ref{NC.eq5}) reduces Eq. (\ref{NC.eq2})
 to the quadratic equation (\textit{characteristic equation})
 \begin{equation}
 \label{NC.eq6}
 \lambda^2 + A \lambda + B = 0.
 \end{equation}

If Eq. (\ref{NC.eq6}) has two distinct roots, $\lambda_1 \not=
\lambda_2,$ then Eq. (\ref{NC.eq5}) provides two linearly
independent solutions
 $$
 y_1 = {\rm e}^{\lambda_1 x}, \quad y_2 = {\rm e}^{\lambda_2 x},
 $$
 and hence, a fundamental set of solutions. If the roots are real,
 the general solution to Eq. (\ref{NC.eq2}) is
 \begin{equation}
 \label{NC.lem1.eq2-sol}
 y(x) = K_1 {\rm e}^{\lambda_1 x} + K_2
 {\rm e}^{\lambda_2 x}, \quad K_1, K_2 = {\rm const.}
 \end{equation}

If the roots are complex, $\lambda_1 = \gamma + i \theta, \
\lambda_2 = \gamma - i \theta,$
 the general solution to Eq. (\ref{NC.eq2}) is given by
 \begin{equation}
 \label{NC.lem1.eq2-compsol}
 y(x) = [K_1 \cos (\theta x) + K_2 \sin (\theta
 x)]\,{\rm e}^{\gamma x}, \quad K_1, K_2 = {\rm const.}
 \end{equation}

 In the case of equal roots $\lambda_1 =\lambda_2,$ standard texts
 in differential equations make a guess, without motivation, that the general solution
 has the form
 \begin{equation}
 \label{NC.lem1.eq2}
 y(x) = (K_1 + K_2\, x)\,{\rm e}^{\lambda_1 x},
 \quad K_1, K_2 = {\rm const.}
 \end{equation}
 The motivation is given in \cite{ibr99}, Section
 13.2.2, and states the following.

 \begin{lem}
 \label{NC.cc.lem1}
 Eq. (\ref{NC.eq2}) can be mapped to the equation $z'' = 0$  by
 a linear change of the dependent variable
 \begin{equation}
 \label{NC.lem1.eq0}
 y = \sigma (x)\, z, \quad \sigma (x) \not= 0,
 \end{equation}
 if and only if the characteristic equation (\ref{NC.eq6}) has equal
 roots. Specifically, if $\lambda_1 =\lambda_2,$
 the substitution
 \begin{equation}
 \label{NC.lem1.eq1}
  y = z \, {\rm e}^{\lambda_1 x}
 \end{equation}
 carries Eq. (\ref{NC.eq2}) to the equation
  $
  z'' = 0.
  $
  Substituting in (\ref{NC.lem1.eq1}) the solution $z = K_1 + K_2\, x$ of the equation
  $z'' = 0,$ we obtain the general solution (\ref{NC.lem1.eq2})
 to Eq. (\ref{NC.eq2}) whose coefficients satisfy the condition of
 equal roots for  Eq. (\ref{NC.eq6}):
 \begin{equation}
 \label{NC.lem1.eq3}
  A^2 - 4 B = 0.
 \end{equation}
 \end{lem}
 {\bf Proof.}
 The reckoning shows that after the substitution (\ref{NC.lem1.eq0})
  Eq. (\ref{NC.eq1}) becomes (see, e.g. \cite{ibr06P}, Section 3.3.2)
  $$
  z'' + I(x) z = 0,
  $$
 where
 $$
 I(x) = g(x) - \frac{1}{4}\,f^2(x) - \frac{1}{2}\, f'(x)
 $$
 and the function $\sigma(x)$ in the transformation
 (\ref{NC.lem1.eq0}) has the form
 \begin{equation}
 \label{NC.lem1.eq4}
 \sigma(x) = {\rm e}^{- \frac{1}{2}\int f(x) dx}.
 \end{equation}
 Hence,  Eq.
 (\ref{NC.eq1}) is carried into equation $z'' = 0$ if and only if
 $I(x) = 0,$ i.e.
 \begin{equation}
 \label{NC.lem1.eq5}
f^2(x) + 2f'(x) - 4 g(x) = 0.
 \end{equation}
 In the case of Eq. (\ref{NC.eq2}), the condition
 (\ref{NC.lem1.eq5}) is identical with Eq. (\ref{NC.lem1.eq3}),
 and the function $\sigma(x)$ given by Eq. (\ref{NC.lem1.eq4})
 becomes
 \begin{equation}
 \label{NC.lem1.eq6}
 \sigma(x) = {\rm e}^{- \frac{A}{2}\,x}.
 \end{equation}
 Furthermore, under the condition Eq. (\ref{NC.lem1.eq3}), the
 repeated root of  Eq. (\ref{NC.eq6}) is $\lambda_1 = - A/2.$
 Therefore Eq. (\ref{NC.lem1.eq6}) can be written
 $$
 \sigma(x) = {\rm e}^{\lambda_1\,x}
 $$
 and we arrive at the substitution (\ref{NC.lem1.eq1}), and hence
 at the solution (\ref{NC.lem1.eq2}), thus proving the lemma.

 \subsection{Euler's equation}
 \label{NC.ee}

 Consider Euler's equation
 \begin{equation}
 \label{NC.eq7}
 x^2\,y'' + A\,  x\, y' + B\, y = 0, \quad A, B = {\rm const.}
 \end{equation}
It is double homogeneous (see \cite{ibr06P}, Section 6.6.1), i.e.
admits the dilation groups in $x$ and in $y$ with the
 generators
 \begin{equation}
 \label{NC.eq8}
 X_1 = x \frac{\partial}{\partial x}\,, \quad
 X_2 = y \frac{\partial}{\partial y}\,\cdot
 \end{equation}
 We proceed as in Section \ref{NC.cc} and find the invariant solutions
 for the linear combination $X = X_1 + \lambda X_2:$
  \begin{equation}
 \label{NC.eq9}
  X = x \frac{\partial}{\partial x} + \lambda y
\frac{\partial}{\partial y}\,, \quad \lambda = {\rm const.}
 \end{equation}
The characteristic equation
 $$
 \frac{d y}{y} = \lambda \, \frac{d x}{x}
 $$
 of the equation $X(J) = 0$ for the invariants $J(x, y)$ yields the
 invariant
 $$
 J = y\,x^{-\lambda}
 $$
 for the operator (\ref{NC.eq9}). The invariant
 solutions are given by $J = C,$ ~whence
 $$
 y = C x^{\lambda}, \quad C = {\rm const.}
 $$
 Due to the homogeneity of Eq. (\ref{NC.eq2}) we can set $C = 1$ and obtain
 \begin{equation}
 \label{NC.eq10}
  y = x^{\lambda}.
 \end{equation}
 Differentiating and multiplying by $x,$ we have:
 $$
  x y' = \lambda\, x^\lambda, \quad  x^2 y'' = \lambda (\lambda- 1) x^\lambda.
 $$
 Substituting in Eq. (\ref{NC.eq9}) and dividing by the common factor $C x^\lambda$
 we obtain the following {\it characteristic equation} for Euler's
 equation (\ref{NC.eq7}):
 \begin{equation}
 \label{NC.eq11}
 \lambda^2 + (A  - 1)\, \lambda + B  = 0.
 \end{equation}

 \begin{rem}
 \label{Eul_cons}
 According to Eqs. (\ref{NC.eq6}), (\ref{NC.eq11}), the characteristic
 equation for Euler's equation written in the form
 \begin{equation}
 \label{NC.eq12}
 x^2\,y'' + (A + 1)\,  x\, y' + B\, y = 0
 \end{equation}
 is identical with the characteristic equation (\ref{NC.eq6})
  for Eq. (\ref{NC.eq2}) with constant coefficients.
  \end{rem}

 \section{New examples of reducible equations}
 \label{NC.so}
 \setcounter{equation}{0}

 \subsection{First example}
 \label{NC.f}

 Let us find the linear second-order
 equations (\ref{NC.eq1}),
 \begin{equation}
 \tag*{(\ref{NC.eq1})}
 y'' + f(x) y' + g(x) y = 0,
 \end{equation}
admitting the operator
 \begin{equation}
 \label{NC.so1}
 X_1 = x^\alpha \frac{\partial}{\partial x}\,,
 \end{equation}
 where $\alpha$ is any real-valued parameter. Taking the second
 prolongation of $X_1,$
 $$
 X_1 = x^\alpha \frac{\partial}{\partial x} -
 \alpha x^{\alpha - 1} y' \frac{\partial}{\partial y'} -
 \left[\alpha(\alpha - 1) x^{\alpha - 2} y' + 2 x^{\alpha - 1}
  y''\right] \frac{\partial}{\partial y''}\,,
 $$
 we write the invariance condition of Eq. (\ref{NC.eq1}),
 $$
 X_1(y'' + f(x) y' + g(x) y)\Big|_{(\ref{NC.eq1})} = 0,
 $$
 and obtain:
 \begin{equation}
 \label{NC.so2}
  x^{\alpha - 2} [x^2 f' + \alpha x f - \alpha(\alpha - 1)] y' +
 x^{\alpha - 1} [x g' + 2 \alpha g]y = 0.
 \end{equation}
 Since Eq. (\ref{NC.so2}) should be satisfied identically in the
 variables $x, y, y',$ it splits into two equations:
 \begin{equation}
 \label{NC.so3}
  x^2 f' + \alpha x f - \alpha(\alpha - 1)= 0, \quad
x g' + 2 \alpha g = 0.
 \end{equation}
Solving the first-order linear differential equations (\ref{NC.so3})
for the unknown functions $f(x)$ and $g(x),$ we obtain:
 $$
 f(x) = \frac{\alpha}{x} + A x^{- \alpha}, \quad
 g(x) = B x^{- 2 \alpha}, \quad A, B = {\rm const.}
 $$
Thus, we arrive at the following linear equation admitting the
operator (\ref{NC.so1}):
 \begin{equation}
 \label{NC.so4}
 x^{2 \alpha} y'' +  \left(A x^\alpha + \alpha x^{2 \alpha -
 1}\right) y' + B y = 0 , \quad A, B = {\rm const.}
 \end{equation}

 \begin{rem}
When $\alpha = 0,$ Eq. (\ref{NC.so4}) yields the equation
(\ref{NC.eq2}) with constant coefficients. When $\alpha = 1,$ Eq.
(\ref{NC.so4}), upon setting $A + 1$ as a new coefficient $A,$
coincides with Euler's equation (\ref{NC.eq7}).
 \end{rem}

Since Eq. (\ref{NC.so4}) is linear homogeneous, it admits, along
with (\ref{NC.so1}), the operator
 $$
  X_2 = y \frac{\partial}{\partial y}\,\cdot
 $$
Now we proceed as in Section \ref{NC.gp} and find the invariant
solutions
 for the linear combination $X = X_1 + \lambda X_2:$
 $$
  X = x^\alpha \frac{\partial}{\partial x} + \lambda y
\frac{\partial}{\partial y}\,, \quad \lambda = {\rm const.}
 $$
Let $\alpha \not= 1.$ The characteristic equation
 $$
 \frac{d y}{y} = \lambda \, \frac{d x}{x^\alpha}
 $$
 of the equation $X(J) = 0$ for the invariants $J(x, y)$ yields the
 invariant
 $$
 J = y\,{\rm e}^{\frac{\lambda}{\alpha - 1}\,x^{1 -\alpha}}.
 $$
  Hence, the invariant
 solutions are obtained by setting $J = C,$ ~whence letting $C =
 1$ we have
 \begin{equation}
 \label{NC.so5}
  y = {\rm e}^{\frac{\lambda}{1 - \alpha}\,x^{1 -\alpha}}.
 \end{equation}
 Differentiating we have:
 $$
  y' = \lambda\, x^{- \alpha}\,{\rm e}^{\frac{\lambda}{1 - \alpha}\,x^{1 -\alpha}},
  \quad  y'' = \big[\lambda^2 x^{- 2\alpha} - \lambda \alpha x^{- \alpha - 1}\big]
   {\rm e}^{\frac{\lambda}{1 - \alpha}\,x^{1 -\alpha}}.
 $$
 Substituting in Eq. (\ref{NC.so4}) and dividing by the non-vanishing factor
 ${\rm e}^{\frac{\lambda}{1 - \alpha}\,x^{1 -\alpha}},$
 we obtain the following {\it characteristic equation} for Eq.
 (\ref{NC.so4}):
 \begin{equation}
 \label{NC.so6}
 \lambda^2 + A\, \lambda + B  = 0.
 \end{equation}
 Eq. (\ref{NC.so6}) is identical with the characteristic equation (\ref{NC.eq6}) for
 the equation (\ref{NC.eq2}) with constant coefficients.

 \begin{rem}
 Eq. (\ref{NC.so5}) contains Euler's substitution (\ref{NC.eq5})
for the equation (\ref{NC.eq2}) with constant coefficients as a
particular case  $\alpha = 0.$
 \end{rem}

 \begin{rem}
Using the statement that  Eq. (\ref{NC.eq1}) is mapped to the
equation $z'' = 0$ if and only if the function
 $$
 I(x) = g(x) - \frac{1}{4}\, f^2(x) -
 \frac{1}{2}\,f'(x)
 $$
 vanishes (see Lemma \ref{NC.cc.lem1}), one can verify that
 Eq. (\ref{NC.so5}) is equivalent by function to the equation $z'' = 0$ if and only if
\begin{equation}
 \label{NC.so7}
 A^2 - 4B = 0 \quad {\rm and} \quad \alpha  = 0 \ {\rm or} \  \alpha  = 2.
 \end{equation}
 The first equation in (\ref{NC.so7}) means that the characteristic equation
 (\ref{NC.so6}) has a repeated root, and hence there is only one
 solution of the form (\ref{NC.so5}). Then, using the reasoning of Lemma \ref{NC.cc.lem1}
  we can show that the general solution to
 Eq. (\ref{NC.so4}) with $A^2 - 4B = 0, \ \alpha = 2,$ i.e. of the equation
 $$
 y'' + \left(\frac{A}{x^2} + \frac{2}{x}\right) y' +
 \frac{B}{x^4}\, y = 0, \quad A^2 - 4B = 0,
 $$
 is given by
 \begin{equation}
 \label{NC.so8}
 y = \left(K_1 + \frac{K_2}{x}\right) {\rm e}^{-
 \frac{\lambda}{x}},
 \end{equation}
 where $K_1, K_2$ are arbitrary constants and $\lambda$ is the
 repeated root of the characteristic equation (\ref{NC.so6}).
 For a more general statement, see Section \ref{NC.gen_exs}.
 \end{rem}

 \subsection{Second example}
 \label{NC.s}

  Let us find the linear second-order
 equations (\ref{NC.eq1}) admitting the projective group with the generator
 \begin{equation}
 \label{NC.so9}
 X_1 = x^2 \frac{\partial}{\partial x} + x y \frac{\partial}{\partial y}\,\cdot
 \end{equation}
 Taking the second prolongation of the operator (\ref{NC.so9}),
 $$
 X_1 = x^2 \frac{\partial}{\partial x} + x y \frac{\partial}{\partial y}
 + (y - x y') \frac{\partial}{\partial y'} -
 3 x y'' \frac{\partial}{\partial y''}\,,
 $$
 and writing the invariance condition of Eq. (\ref{NC.eq1}),
 $$
 X_1(y'' + f(x) y' + g(x) y)\Big|_{(\ref{NC.eq1})} = 0,
 $$
 we obtain:
 \begin{equation}
 \label{NC.so10}
  x(x f' + 2 f) y' + (x^2 g' + 4 x g + f)y = 0.
 \end{equation}
 Eqs. (\ref{NC.so10}) yield:
 $$
 f(x) = \frac{A}{x^2}\,, \quad
 g(x) = \frac{B}{x^4} - \frac{A}{x^3}\,, \quad A, B = {\rm const.}
 $$
Thus, we arrive at the following linear equation admitting the
operator (\ref{NC.so9}):
 \begin{equation}
 \label{NC.so11}
 x^4 y'' + A x^2 y' + (B - Ax) y = 0 , \quad A, B = {\rm const.}
 \end{equation}

Eq. (\ref{NC.so11}) is linear homogeneous, and hence admits, along
with (\ref{NC.so9}), the operator
 $$
  X_2 = y \frac{\partial}{\partial y}\,\cdot
 $$
Now we proceed as in Section \ref{NC.gp} and find the invariant
solutions
 for the linear combination $X = X_1 + \lambda X_2:$
 $$
  X = x^2 \frac{\partial}{\partial x} + (x + \lambda) y
 \frac{\partial}{\partial y}\,, \quad \lambda = {\rm const.}
 $$
 Rewriting the characteristic equation of the equation $X(J) = 0$ for the invariants $J(x, y)$ in the form
 $$
 \frac{d y}{y} =  \frac{x + \lambda}{x^2} d x
 $$
we obtain the invariant
 $$
 J = \frac{y}{x}\,{\rm e}^{\frac{\lambda}{x}}.
 $$
 Setting $J = C$ and letting $C = 1$ we obtain the following form of the invariant
 solutions:
 \begin{equation}
 \label{NC.so12}
  y = x\,{\rm e}^{- \frac{\lambda}{x}}.
 \end{equation}
Substituting (\ref{NC.so12}) in Eq. (\ref{NC.so11}) we reduce the
differential equation (\ref{NC.so11}) to the algebraic equation
 $$
 \lambda^2 + A\, \lambda + B  = 0
 $$
 which is identical with the characteristic equation (\ref{NC.eq6}) for
 the equation (\ref{NC.eq2}) with constant coefficients.

 \begin{exa}
 Solve the equation
 \begin{equation}
 \label{NC.so13}
 y'' + \frac{\omega^2}{x^4}\,y = 0, \quad \omega = {\rm const.}
 \end{equation}
 This is an equation of the form (\ref{NC.so11}) with $A = 0, \ B
 = \omega^2.$ The algebraic equation (\ref{NC.eq6}) yields
 $\lambda_1 = - i \omega, \ \lambda_2 = i \omega,$ and hence we
 have two independent invariant solutions (\ref{NC.so12}):
 $$
  y_1 = x\,{\rm e}^{i \frac{\omega}{x}}, \quad
 y_2 = x\,{\rm e}^{- i \frac{\omega}{x}}.
 $$
 Taking their real and imaginary parts, just like in the case of
 constant coefficient equations, we obtain the following
 fundamental system of solutions:
 \begin{equation}
 \label{NC.so14}
 y_1 = x\,\cos \left(\frac{\omega}{x}\right), \quad
 y_2 = x\,\sin \left(\frac{\omega}{x}\right).
 \end{equation}
 Hence, the general solution to Eq. (\ref{NC.so12}) is given by
 $$
 y = x\,\Big[C_1 \cos \left(\frac{\omega}{x}\right)
 + C_2 \sin \left(\frac{\omega}{x}\right)\Big].
 $$
 \end{exa}

We can also solve the non-homogeneous equation
 \begin{equation}
 \label{NC.so15}
 y'' + \frac{\omega^2}{x^4}\,y = F(x),
 \end{equation}
 e.g. by the method of \textit{variation of parameters}, and
 obtain
 \begin{align}
 \label{NC.so16}
  & y = x\,\Big[C_1 \cos \left(\frac{\omega}{x}\right) + C_2 \sin \left(\frac{\omega}{x}\right)\Big]\\[2ex]
  & +
 \frac{x}{\omega}\, \Big[\cos \left(\frac{\omega}{x}\right)
  \int x F(x) \sin \left(\frac{\omega}{x}\right) dx
 - \sin \left(\frac{\omega}{x}\right) \int x F(x) \cos \left(\frac{\omega}{x}\right) dx \Big].\notag
 \end{align}

 \section{General result for second-order equations}
 \label{NC.gen}
 \setcounter{equation}{0}

 \subsection{Main statements}

 Note that the operators given in Eqs. (\ref{NC.eq3}),
 (\ref{NC.eq8}), (\ref{NC.so1}), (\ref{NC.so9}) are particular
 cases of the generator
 \begin{equation}
 \label{NC.gen1}
 X_1 = \phi(x) \frac{\partial}{\partial x} + \sigma(x) y \frac{\partial}{\partial y}
 \end{equation}
 of the general equivalence group of all linear ordinary
 differential equations. We will find now all linear second-order
 equations (\ref{NC.eq1}) admitting the operator (\ref{NC.gen1})
 with any \textit{fixed} functions $\phi(x)$ and $\sigma(x).$

 Taking the second prolongation of the operator (\ref{NC.gen1}),
 $$
 X_1 = \phi \frac{\partial}{\partial x} + \sigma y \frac{\partial}{\partial y}
 + \big[\sigma ' y + (\sigma - \phi') y'\big] \frac{\partial}{\partial y'}
 + \big[\sigma '' y + (2 \sigma' - \phi'') y' + (\sigma - 2 \phi') y'' \big]
  \frac{\partial}{\partial y''}\,,
 $$
 and writing the invariance condition of Eq. (\ref{NC.eq1}),
 $$
 X_1(y'' + f(x) y' + g(x) y)\Big|_{(\ref{NC.eq1})} = 0,
 $$
 we obtain:
 $$
  (\phi f' + f \phi' + 2 \sigma' - \phi'') y' + (\phi g' + 2 \phi' g +
   f \sigma' + \sigma'')y = 0.
 $$
 It follows:
 \begin{align}
 \label{NC.gen2}
 & \phi f' + f \phi' \ + 2 \sigma' - \phi''= 0, \notag\\
 & \phi g' + 2 \phi' g + f \sigma' + \sigma'' = 0.
 \end{align}
 The first equation (\ref{NC.gen2}) is written
 $$
 (\phi f)' = (\phi' - 2 \sigma)'
 $$
 and yields:
  \begin{equation}
 \label{NC.gen3}
 f(x) = \frac{1}{\phi}\,\big[A + \phi' - 2 \sigma\big], \quad A = {\rm
 const.}
 \end{equation}
 Substituting this in the second equation (\ref{NC.gen2}), we
 obtain the following non-homogeneous linear first-order equation for determining
 $g(x):$
 \begin{equation}
 \label{NC.gen4}
 \phi g' + 2 \phi' g = - \sigma'' - \frac{\sigma'}{\phi}\,
 \big[A + \phi' - 2 \sigma\big].
 \end{equation}
 The homogeneous equation
 $$
 \phi(x) g' + 2 \phi'(x) g = 0
 $$
 with a given function $\phi(x)$ yields
 $$
 g = \frac{C}{\phi^2(x)}\,\cdot
 $$
By variation of the parameter $C,$ we set
 $$
 g = \frac{u(x)}{\phi^2(x)}\,,
 $$
 substitute it in Eq. (\ref{NC.gen4}) and obtain:
 $$
 u' = - A \sigma' + 2 \sigma \sigma' - \phi' \sigma' - \phi
 \sigma'' \equiv - (A \sigma)' + (\sigma^2)' - (\phi \sigma')',
 $$
whence
 $$
 u = B - A \sigma + \sigma^2 - \phi \sigma', \quad B = {\rm const.}
 $$
 Therefore,
 \begin{equation}
 \label{NC.gen5}
 g = \frac{1}{\phi^2(x)}\, \Big[B - A \sigma + \sigma^2 - \phi \sigma'\Big].
 \end{equation}
Thus, we have arrived at the following result.

 \begin{thm}
 \label{geu.t1}
 The  homogeneous linear second-order
 equations (\ref{NC.eq1}) admitting the operator (\ref{NC.gen1})
 with any given functions $\phi = \phi(x)$ and $\sigma = \sigma(x)$
 have the form
 \begin{equation}
 \label{NC.gen6}
 \phi^2\,y'' + \big(A + \phi' - 2 \sigma\big) \phi\, y' +
  \big(B - A \sigma + \sigma^2 - \phi \sigma'\big) y = 0.
 \end{equation}
 \end{thm}

 Now we use the homogeneity of Eq. (\ref{NC.gen6})
characterized by the generator
 $$
  X_2 = y \frac{\partial}{\partial y}\,\cdot
 $$
Namely, we look for the invariant solutions with respect to the
linear combination $X = X_1 + \lambda X_2:$
 $$
  X = \phi(x) \frac{\partial}{\partial x} + (\sigma(x) + \lambda) y
 \frac{\partial}{\partial y}\,, \quad \lambda = {\rm const.,}
 $$
 and arrive at the following statement reducing the problem of integration of the
 differential equation (\ref{NC.gen6}) to solution of the quadratic equation,
 namely, the characteristic equation as in the case of equations
 with constant coefficients.

 \begin{thm}
  \label{geu.t2}
  Eq. (\ref{NC.gen6}) has the invariant solutions of the form
 \begin{equation}
 \label{NC.gen7}
 y = {\rm e}^{\int \frac{\sigma (x) + \lambda}{\phi(x)}\,dx},
 \end{equation}
 where $\lambda$ satisfies the \textit{characteristic equation}
 \begin{equation}
 \label{NC.gen8}
 \lambda^2 + A \lambda + B = 0.
 \end{equation}
 \end{thm}
 {\bf Proof.} We solve the equation
   $X(J) = 0$ for the invariants $J(x, y),$ i.e. integrate the
   equation
 $$
 \frac{d y}{y} =  \frac{\sigma (x) + \lambda}{\phi(x)} d x
 $$
and obtain the invariant
 $$
 J = y\,{\rm e}^{- \int \frac{\sigma (x) + \lambda}{\phi(x)}\,dx}.
 $$
 Setting $J = C$ and letting $C = 1$ we obtain Eq. (\ref{NC.gen7}) for the invariant
 solutions. Thus, we have:
 \begin{align}
 \label{NC.gen9}
 & y = {\rm e}^{\int \frac{\sigma (x) + \lambda}{\phi(x)}\,dx}, \quad
 y' = \frac{\sigma  + \lambda}{\phi}\,{\rm e}^{\int \frac{\sigma (x) + \lambda}{\phi(x)}\,dx}\\[1.5ex]
 & y'' = \frac{1}{\phi^2}\,\left[(\sigma + \lambda)^2 - (\sigma + \lambda)\phi'
 + \phi \sigma'\,\right]\,{\rm e}^{\int \frac{\sigma (x) + \lambda}{\phi(x)}\,dx}. \notag
 \end{align}
Substituting (\ref{NC.gen9}) in  Eq. (\ref{NC.gen6}) we obtain Eq.
(\ref{NC.gen8}), thus completing the proof.

 \subsection{Distinct roots of the characteristic equation}

 It is manifest that if the characteristic equation (\ref{NC.gen8}) has
 distinct real roots $\lambda_1 \not= \lambda_2,$ the general solution
 to  Eq. (\ref{NC.gen6}) is given by
 \begin{equation}
 \label{NC.gen6_sol}
 y(x) = K_1\,{\rm e}^{\int \frac{\sigma (x) + \lambda_1}{\phi(x)}\,dx}
  + K_2\,{\rm e}^{\int \frac{\sigma (x) +
  \lambda_2}{\phi(x)}\,dx}, \quad K_1, K_2 = {\rm const.}
 \end{equation}

 In the case of complex roots, $\lambda_1 = \gamma + i \theta, \
 \lambda_2 = \gamma - i \theta,$
 the general solution to  Eq. (\ref{NC.gen6}) is given by
 \begin{equation}
 \label{NC.gen6_compsol}
 y(x) = \left[K_1\,\cos \left(\theta \int \frac{dx}{\phi (x)}\right)
 + K_2\,\sin \left(\theta \int \frac{dx}{\phi (x)}\right)\right]
  {\rm e}^{\int \frac{\sigma (x) + \gamma}{\phi(x)}\,dx}.
 \end{equation}

 \subsection{The case of repeated roots}

 \begin{thm}
 \label{geu.t3}
 If the characteristic equation (\ref{NC.gen8}) has equal roots
 $\lambda_1 = \lambda_2,$  the general solution to  Eq. (\ref{NC.gen6})
 is given by
 \begin{equation}
 \label{NC.gen10}
 y = \left[K_1 + K_2 \, \int \frac{dx}{\phi(x)}\right]
 {\rm e}^{\int \frac{\sigma (x) + \lambda_1}{\phi(x)}\,dx}, \quad
 K_1, K_2 = {\rm const.}
 \end{equation}
 \end{thm}
{\bf Proof.} Let us  new variables $t$ and  $z$ defined by the
linear first-order equations
 \begin{equation}
 \label{NC.gen11}
 X_1(t) \equiv \phi(x) \frac{\partial t}{\partial x} +
 \sigma(x)\, y\, \frac{\partial t}{\partial y} = 1, \quad
 X_2(t) \equiv y\, \frac{\partial t}{\partial y} = 0
 \end{equation}
 and
 \begin{equation}
 \label{NC.gen12}
 X_1(z) \equiv \phi(x)\, \frac{\partial z}{\partial x} +
 \sigma(x)\, y\, \frac{\partial z}{\partial y} = 0, \quad
 X_2(z) \equiv y\, \frac{\partial z}{\partial y} = z,
 \end{equation}
 respectively. Eqs. (\ref{NC.gen11}) are easily solved and yield
 \begin{equation}
 \label{NC.gen13}
t = \int \frac{dx}{\phi(x)}\,\cdot
 \end{equation}
Integration of the second equation (\ref{NC.gen12}) with respect to
$y$ gives
 $$
 z = v(x) y.
 $$
Substituting this in the first equation (\ref{NC.gen12}) we obtain
 $$
 \phi(x)\, \frac{dv}{d x} + \sigma(x)\,v = 0, \quad {\rm whence}
 \quad v = {\rm e}^{- \int \frac{\sigma(x)}{\phi(x)}\,d x}.
 $$
 %whence
 %$$
 %v(x) = {\rm e}^{- \int \frac{\sigma(x)}{\phi(x)}\,d x}.
 %$$
Thus,
 \begin{equation}
 \label{NC.gen14}
 z = y \, {\rm e}^{- \int \frac{\sigma(x)}{\phi(x)}\,d x}.
 \end{equation}

The passage to the new variables (\ref{NC.gen13}), (\ref{NC.gen14})
converts the operator $X_1$ given by (\ref{NC.gen1}) to the
translation generator without changing the form of the dilation
generator $X_2.$ In other words, upon introducing the new
independent and dependent variables $t$ and $z$ given by
(\ref{NC.gen13}) and (\ref{NC.gen14}), respectively, we arrive at
the operators (\ref{NC.eq3}). Hence, in the new variables, Eq.
(\ref{NC.gen14}) becomes an equation with constant coefficients.
Invoking that the equations (\ref{NC.gen6}) and (\ref{NC.eq2}) have
Eq. (\ref{NC.gen8}) as their common characteristic equation, we use
Lemma \ref{NC.cc.lem1} and write
$$
z = (K_1 + K_2\, t)\,{\rm e}^{\lambda_1 t}.
$$
Substituting this in Eq. (\ref{NC.gen14}) and replacing $t$ and $z$
by their expressions (\ref{NC.gen13}) and (\ref{NC.gen14}),
respectively, and solving for $y,$ we finally arrive at Eq.
(\ref{NC.gen10}).

 \begin{rem}
 \label{gen_nh}
We can easily solve the non-homogeneous equation Eq.
(\ref{NC.gen6}):
 \begin{equation}
 \label{NC.gen6_nh}
 \phi^2\,y'' + \big(A + \phi' - 2 \sigma\big) \phi\, y' +
  \big(B - A \sigma + \sigma^2 - \phi \sigma'\big) y = F(x).
 \end{equation}
 Namely, we rewrite Eq. (\ref{NC.gen6}) in the form
   $$
  y'' + a(x)\, y' + b(x)\, y = P(x)
  $$
and employ the representation of the general solution (see, e.g.
\cite{ibr06P}, Section 3.3.4)
 \begin{equation}
 \label{NC.gen6_nhsol}
 y = K_1 y_1(x) + K_2 y_2(x) - y_1(x)\int \frac{y_2 (x) P(x)}{W(x)}\,d x
 + y_2(x)\int \frac{y_1 (x) P(x)}{W(x)}\,d x
  \end{equation}
  furnished by the method of variation of parameters. Here
  $$
  W(x) = y_1(x) y_2'(x) - y_2(x) y_1'(x)
  $$
  is the Wronskian of a fundamental
 system of solutions $y_1(x), ~y_2(x)$ for the homogeneous equation
 $$
 y'' + a(x)\, y' + b(x)\, y = 0.
 $$
 \end{rem}

 \section{Examples to Section \ref{NC.gen}}
  \label{NC.gen_exs}
 \setcounter{equation}{0}

Euler's substitution (\ref{NC.eq5}) as well as the solutions
(\ref{NC.eq10}), (\ref{NC.so5}) and (\ref{NC.so12}) are encapsulated
in Eq. (\ref{NC.gen7}). We will consider now these and several other
examples.

 \begin{exa}
 \label{NC.gen_ex1}
 Let us take $\phi (x) = 1, \ \sigma (x) = 0.$ Then Eqs. (\ref{NC.gen6}), (\ref{NC.gen7})
 and (\ref{NC.gen10}) coincide with Eqs. (\ref{NC.eq2}), (\ref{NC.eq5})
 and (\ref{NC.lem1.eq2}), respectively. Eq. (\ref{NC.gen6_compsol}) becomes
  (\ref{NC.lem1.eq2-compsol}).
 \end{exa}

 \begin{exa}
 \label{NC.gen_ex2}
 Let us take $\phi(x) = x, \ \sigma (x)= 0.$ Then Eq. (\ref{NC.gen6}) becomes Euler's
 equation written in the form (\ref{NC.eq12}), Eq. (\ref{NC.gen7}) yields Eq.
 (\ref{NC.eq10}) for invariant solutions, whereas Eq. (\ref{NC.gen10}) provides the general solution
 \begin{equation}
 \label{NC.gen_exs.eq1}
 y(x) = (K_1 + K_2\, \ln x)\,x^{\lambda_1}
 \quad K_1, K_2 = {\rm const.,}
 \end{equation}
 to Euler's  equation (\ref{NC.eq12}) whose characteristic  equation (\ref{NC.eq6})
has equal roots. Eq. (\ref{NC.gen6_compsol}) leads
 to the following solution for complex roots $\lambda_1 = \gamma + i \theta,
 \lambda_2 = \gamma - i \theta:$
\begin{equation}
 \label{NC.gen_exs.eq1comp}
 y(x) = [K_1 \cos(\theta \ln x) + K_2 \sin(\theta \ln
 x)]\,x^{\gamma}.
 \end{equation}
 \end{exa}

 \begin{exa}
 \label{NC.gen_ex3}
 Let us take $\phi(x) = x^\alpha, \ \sigma (x)= 0.$ Then Eqs. (\ref{NC.gen6}) and
 (\ref{NC.gen7}) coincide with Eqs. (\ref{NC.so4}) and (\ref{NC.so5}),
 respectively, whereas Eq. (\ref{NC.gen10}) provides the
 general solution
 \begin{equation}
 \label{NC.so8gen}
 y(x) = \big(K_1 + K_2\, x^{1 - \alpha}\big){\rm e}^{\frac{\lambda_1}{1 - \alpha}\,x^{1
 -\alpha}},
 \quad K_1, K_2 = {\rm const.,}
 \end{equation}
 to Eq. (\ref{NC.so4}) whose characteristic
equation (\ref{NC.so6}) has  equal roots $\lambda_1 = \lambda_2.$
 Eq. (\ref{NC.so8gen})  extends
 the solution (\ref{NC.so8}) to all equations (\ref{NC.so4}) with with the coefficients
 $A, B, C$ satisfying the condition (\ref{NC.lem1.eq3}) of equal roots
 for the characteristic  equation (\ref{NC.so6}).
 \end{exa}

 \begin{exa}
 \label{NC.gen_ex4}
 Let us take $\phi(x) = 1 + x^2, \ \sigma (x)= x.$ Then Eq. (\ref{NC.gen6}) becomes
 \begin{equation}
 \label{NC.gen_ex4.eq1}
 \big(1 + x^2\big)^2\, y'' +  \big(1 + x^2\big) A\, y'
 + (B - A x - 1)y = 0.
 \end{equation}
  Working out the integral in Eq. (\ref{NC.gen7}),
$$
 \int \frac{\sigma (x) + \lambda}{\phi(x)}\,dx =
 \int \frac{x}{1 + x^2}\,dx + \int \frac{\lambda}{1 + x^2}\,dx =
 \ln \sqrt{1 + x^2}\, + \lambda\, \arctan x,
$$
we obtain the following expression for  the invariant solutions:
 \begin{equation}
 \label{NC.gen_ex4.eq2}
 y = \sqrt{1 + x^2} \ {\rm e}^{\lambda\, \arctan x},
 \end{equation}
 where $\lambda$ satisfies the characteristic equation
 (\ref{NC.gen8}):
 \begin{equation}
 \label{NC.gen_ex4.eq3}
 \lambda^2 + A \lambda + B = 0.
 \end{equation}

 If the characteristic equation (\ref{NC.gen_ex4.eq3}) has distinct
 real roots, $\lambda_1 \not= \lambda_2,$
 the general solution to  Eq. (\ref{NC.gen_ex4.eq1}) is given by
 \begin{equation}
 \label{NC.gen_ex4.eq4}
 y(x) = \sqrt{1 + x^2} \ \left[K_1\,{\rm e}^{\lambda_1\, \arctan x}
 + K_2\, {\rm e}^{\lambda_2\, \arctan x}\right].
 \end{equation}

 In the case of complex roots, $\lambda_1 = \gamma + i \theta, \
 \lambda_2 = \gamma - i \theta,$
 the general solution to  Eq. (\ref{NC.gen_ex4.eq1}) is given by
 \begin{equation}
 \label{NC.gen_ex4.eq5}
 y(x) = \left[K_1\,\cos \left(\theta \, \arctan x\right)
 + K_2\,\sin \left(\theta \, \arctan x\right)\right]
 \sqrt{1 + x^2} \ {\rm e}^{\lambda\, \arctan x}.
 \end{equation}

Finally, if the characteristic equation (\ref{NC.gen_ex4.eq3}) has
equal roots
 $\lambda_1 = \lambda_2,$  the general solution to  Eq. (\ref{NC.gen_ex4.eq1})
 is given by
 \begin{equation}
 \label{NC.gen_ex4.eq6}
 y(x) = \big(K_1 + K_2\, \arctan x\big)\sqrt{1 + x^2} \
  {\rm e}^{\lambda\, \arctan x}.
 \end{equation}
 \end{exa}

 \begin{exa}
 \label{NC.gen_ex5}
 Consider Eq. (\ref{NC.gen_ex4.eq1}) with $A = 0, \ B = \omega^2.$
 Then, according to Example \ref{NC.gen_ex4}, the equation
 \begin{equation}
 \label{NC.gen_ex4.eq7}
 \big(1 + x^2\big)^2\, y'' + (\omega^2 - 1) y = 0
 \end{equation}
 has the following general solution:
 \begin{equation}
 \label{NC.gen_ex4.eq8}
 y(x) = \left[K_1\,\cos \left(\omega \, \arctan x\right)
 + K_2\,\sin \left(\omega \, \arctan x\right)\right]
 \sqrt{1 + x^2}\,.
 \end{equation}
 \end{exa}

 \begin{exa}
 \label{NC.gen_ex6}
 Let us solve the non-homogeneous equation
 \begin{equation}
 \label{NC.gen_ex4.eq9}
 \big(1 + x^2\big)^2\, y'' + (\omega^2 - 1) y = F(x).
 \end{equation}
Example \ref{NC.gen_ex5} provides the fundamental system of
solutions
 $$
 y_1 = \sqrt{1 + x^2}\, \cos \left(\omega \, \arctan x\right),
 \quad y_2 = \sqrt{1 + x^2}\,
 \sin \left(\omega \, \arctan x\right)
 $$
  for the homogeneous equation (\ref{NC.gen_ex4.eq7}). We have:
 $$
 y_1' = \frac{1}{\sqrt{1 + x^2}}\,\left[x\,\cos \left(\omega \, \arctan x\right)
 - \omega \sin \left(\omega \, \arctan x\right)\right],
 $$
 $$
 y_2' = \frac{1}{\sqrt{1 + x^2}}\,\left[x\,\sin \left(\omega \, \arctan x\right)
 + \omega \cos \left(\omega \, \arctan x\right)\right].
 $$
 Hence the Wronskian is $W[y_1, y_2] = y_1 y_2' - y_2 y_1' = \omega.$
 Now we rewrite Eq. (\ref{NC.gen_ex4.eq9}), in accordance with Remark \ref{gen_nh}, in the form
  \begin{equation}
 \label{NC.gen_ex4.eq10}
 y'' + \frac{\omega^2 - 1}{(1 + x^2)^2}\,  y =
 \frac{F(x)}{(1 + x^2)^2}\,,
 \end{equation}
 use Eq. (\ref{NC.gen6_nhsol}) and obtain the following general solution to Eq. (\ref{NC.gen_ex4.eq9}):
 \begin{align}
 \label{NC.gen_ex4.eq11}
 & y(x) = \sqrt{1 + x^2}\,\Big[K_1\,\cos \left(\omega \, \arctan x\right)
 + K_2\,\sin \left(\omega \, \arctan x\right)\notag\\[1ex]
 & - \frac{1}{\omega} \, \cos \left(\omega \, \arctan x\right)
 \int \frac{F(x)}{(1 + x^2)^{3/2}}\, \sin \left(\omega \, \arctan
 x\right)dx\\[1ex]
 & + \frac{1}{\omega}\, \sin \left(\omega \, \arctan x\right)
 \int \frac{F(x)}{(1 + x^2)^{3/2}}\, \cos \left(\omega \, \arctan x\right) dx \Big].\notag
 \end{align}
 \end{exa}

 \section{Third-order equations}
  \label{NC.gen_third}
 \setcounter{equation}{0}

 The previous results can be extended to higher-order linear ordinary differential
  equations. I will discuss here the third-order equations
  \begin{equation}
  \label{NC.third_gen1}
 y''' + f(x) y'' + g(x) y' + h(x) y = 0.
 \end{equation}

 \begin{thm}
 \label{geu.t4}
 The homogeneous linear third-order
 equations (\ref{NC.third_gen1}) admitting the operator (\ref{NC.gen1})
 with any given  $\phi = \phi(x)$ and $\sigma = \sigma(x)$ have
 the form
 \begin{align}
 \label{NC.third_gen2}
 & \phi^3\,y''' + \big[A + 3 (\phi' - \sigma)\big] \phi^2\, y''\\[1ex]
 & + \big[B + A\phi' - 2 A \sigma + \phi \phi'' + (\phi')^2 - 3 (\phi \sigma)'
   + 3 \sigma^2\big] \phi\, y'\notag\\[1ex]
 &  +  \big[C - B \sigma + A \sigma^2 -  A \phi \sigma' - \sigma^3
   - \phi^2 \sigma '' - \phi \phi' \sigma' + 3 \phi \sigma \sigma'
   \big]\, y = 0.\notag
 \end{align}
 \end{thm}
 {\bf Proof.} We take the third prolongation of the operator
 (\ref{NC.gen1}):
 \begin{align}
  X_1 = & \ \phi \frac{\partial}{\partial x} + \sigma y \frac{\partial}{\partial y}
 + \big[\sigma ' y + (\sigma - \phi') y'\big] \frac{\partial}{\partial y'}\notag\\[1.5ex]
  & + \big[\sigma '' y + (2 \sigma' - \phi'') y' + (\sigma - 2 \phi') y'' \big]
  \frac{\partial}{\partial y''}\notag\\[1.5ex]
 & + \big[\sigma''' y + (3 \sigma'' - \phi''') y' + 3 (\sigma' - \phi'') y''
 + (\sigma - 3 \phi') y''' \big]
  \frac{\partial}{\partial y'''}\,,\notag
 \end{align}
and write the invariance condition of Eq. (\ref{NC.third_gen1}):
 $$
 X_1(y''' + f(x) y'' + g(x) y' + h(x) y)\Big|_{(\ref{NC.third_gen1})} = 0.
 $$
 We annul the coefficients for $y'', y'$ and $y$ of the left-hand side
 of the above equation and split it into the following three equations:
  \begin{align}
 & \phi f' + \phi' f + 3 (\sigma' - \phi'') = 0, \label{NC.third_gen3}\\[1ex]
 & \phi g' + 2 \phi' g  + (2 \sigma' - \phi'') f - \phi'''
  + 3 \sigma'' = 0,\label{NC.third_gen4}\\[1ex]
 & \phi h' + 3 \phi' h  + \sigma'' f + \sigma' g + \sigma''' = 0.\label{NC.third_gen5}
 \end{align}
Eq. (\ref{NC.third_gen3}) is written
 $$
 (\phi f)' = 3 (\phi' - \sigma)'
 $$
 and yields:
 \begin{equation}
 \label{NC.third_gen6}
 f = \frac{1}{\phi}\,\big[A + 3(\phi' - \sigma)\big], \quad A = {\rm
 const.}
 \end{equation}
 We substitute (\ref{NC.third_gen6}) in Eq.
 (\ref{NC.third_gen4}), integrate the resulting non-homogeneous
 linear first-order equation for $g$ and obtain:
 \begin{equation}
 \label{NC.third_gen7}
 g = \frac{1}{\phi^2}\,\big[B + A\phi' - 2 A \sigma + \phi \phi'' + (\phi')^2 - 3 (\phi \sigma)'
   + 3 \sigma^2\big],
 \end{equation}
where $B$ is an arbitrary constant. Now we substitute
(\ref{NC.third_gen6}), (\ref{NC.third_gen7}) in Eq.
 (\ref{NC.third_gen5}), integrate the resulting  first-order equation for $h$ and obtain:
 \begin{equation}
 \label{NC.third_gen8}
 h = \frac{1}{\phi^3}\,\big[C - B \sigma + A \sigma^2 -  A \phi \sigma' - \sigma^3
   - \phi^2 \sigma '' - \phi \phi' \sigma' + 3 \phi \sigma \sigma'\big],
 \end{equation}
where $C$ is an arbitrary constant. Finally, substituting
(\ref{NC.third_gen6}), (\ref{NC.third_gen7}) and
(\ref{NC.third_gen8}) in Eq. (\ref{NC.third_gen1}), we arrive at Eq.
(\ref{NC.third_gen2}).

 \begin{thm}
  \label{geu.t5}
  Eq. (\ref{NC.third_gen2}) has the invariant solutions of the
  form (\ref{NC.gen7}),
 \begin{equation}
 \tag*{(\ref{NC.gen7})}
 y = {\rm e}^{\int \frac{\sigma (x) + \lambda}{\phi(x)}\,dx},
 \end{equation}
 where $\lambda$ satisfies the \textit{characteristic equation}
\begin{equation}
 \label{NC.third_gen9}
 \lambda^3 + A \lambda^2 + B \lambda + C = 0.
 \end{equation}
 \end{thm}
 {\bf Proof.} Adding to Eqs. (\ref{NC.gen9}) the expression for the third derivative
 $y'''$ and substituting in  Eq. (\ref{NC.third_gen2}) we obtain Eq.
(\ref{NC.third_gen9}).

 \chapter{Using connection between Riccati and
 second-order linear equations}
 %\chapter{Using connection of Riccati equations with
 %second-order linear equations}
  \label{Ricc.Lin}
    \setcounter{section}{10}
% \addcontentsline{toc}{chapter}{Riccati equation}

 \section{Introduction}
% \addcontentsline{toc}{section}{Introduction}
 \label{Ricc.Lin_int}
 \setcounter{equation}{0}

 Recall that the Riccati equation
 \begin{equation}
 \label{Ricc.Lin_eq1}
 y' = P(x) + Q(x) y + R(x) y^2, \quad R(x)\not= 0,
 \end{equation}
 is mapped by the substitution
 \begin{equation}
 \label{Ricc.Lin_eq2}
 y = - \frac{1}{R(x)}\,\frac{u'}{u}
 \end{equation}
to the linear second-order equation
 \begin{equation}
 \label{Ricc.Lin_eq3}
 u'' + f(x)u' + g(x) u = 0
 \end{equation}
 with the coefficients
 \begin{equation}
 \label{Ricc.Lin_eq4}
 f(x) = - \bigg[Q(x) + \frac{R'(x)}{R(x)}\bigg]\,, \quad   g(x)  =  P(x) R (x).
 \end{equation}
Indeed, we have:
 $$
 y' = - \frac{1}{R}\,\frac{u''}{u} + \frac{R'}{R^2}\,\frac{u'}{u} +
  \frac{1}{R}\,\frac{u'^2}{u^2}
 $$
 and
 $$
 P + Q y + R y^2 = P - \frac{Q}{R}\,\frac{u'}{u} +
 \frac{1}{R}\,\frac{u'^2}{u^2}\,\cdot
 $$
 Substituting these expressions in Eq. (\ref{Ricc.Lin_eq1}) and
 multiplying by $- R u$ we obtain the equation
 $$
 u'' - \frac{R'}{R}\,u' = Q u' - P R u,
 $$
 i.e. Eq. (\ref{Ricc.Lin_eq3}) with the coefficients (\ref{Ricc.Lin_eq4}).

 \section{From Riccati to second-order equations}
 %\section{Integration of second-order linear equations using linearizable Riccati equations}
 \label{Ricc.Lin_sec1}
 \setcounter{equation}{0}

 \subsection{Application to  Equation (\ref{gai.ric2})}
 \label{Ricc.Lin_subsec1}

 Applying Eqs. (\ref{Ricc.Lin_eq3})-(\ref{Ricc.Lin_eq4})
 to Eq. (\ref{gai.ric2}),
 \begin{equation}
 \tag*{(\ref{gai.ric2})}
 y' = Q(x)y + R(x)y^2,
 \end{equation}
 we obtain the following second-order linear equation:
 \begin{equation}
 \label{Ricc.Lin_subsec1.eq1}
 u'' = \bigg(Q + \frac{R'}{R}\bigg)u'.
 \end{equation}
 The integration yields:
 $$
 \ln u' = \int \bigg(Q + \frac{R'}{R}\bigg)dx +  \ln C_1 = \int Q dx + \ln R +  \ln
 C_1.
 $$
 Hence
 \begin{equation}
 \label{Ricc.Lin_subsec1.eq2}
 u' = C_1 R(x) {\rm e}^{\int Q(x)dx}
 \end{equation}
 and
 \begin{equation}
 \label{Ricc.Lin_subsec1.eq3}
 u = C_1 \int R(x) {\rm e}^{\int Q(x)dx}\,dx + C_2.
 \end{equation}
 Substituting (\ref{Ricc.Lin_subsec1.eq2}) and (\ref{Ricc.Lin_subsec1.eq3})
 in Eq. (\ref{Ricc.Lin_eq2}) and denoting $C = C_2/ C_1,$ we arrive
 at the solution (\ref{gai.eq10}) to Eq. (\ref{gai.ric2}):
 \begin{equation}
 \tag*{(\ref{gai.eq10})}
   y = - \frac{{\rm e}^{\int Q(x) dx}}{C + \int R(x) {\rm e}^{\int Q(x)
   dx}\,dx}\,\cdot
 \end{equation}

 \subsection{Application to  Equation (\ref{gai.ric3})}
 \label{Ricc.Lin_subsec2}

The following examples clarify how to use the linearizable equations
(\ref{gai.ric3}),
 \begin{equation}
 \tag*{(\ref{gai.ric3})}
 y' = P(x) + Q(x)y + k[Q(x)- kP(x)]y^2, \quad k = {\rm const.,}
 \end{equation}
 for integrating the corresponding
second-order linear equations (\ref{Ricc.Lin_eq3}).

 \begin{exa}
 \label{Ricc.Lin_sec1.ex1}
 If we apply Eqs. (\ref{Ricc.Lin_eq3})-(\ref{Ricc.Lin_eq4}) to Eq. (\ref{gai.eq17})
 from Example \ref{gai.sub2.ex1},
 $$
 y' = x + 2 xy + x y^2,
 $$
 we obtain the following second-order
 linear equation:
 \begin{equation}
 \label{Ricc.Lin_sec1.eq1}
 u'' - \bigg(2x + \frac{1}{x}\bigg) u' + x^2 u = 0.
 \end{equation}
 Let us integrate this equation. Writing Eq. (\ref{Ricc.Lin_eq2}) in the form
 $$
 \frac{u'}{u} = - R(x) y,
 $$
 substituting here $R(x) = x$  and the expression (\ref{gai.eq17sol}) for $y,$ we
 get
 $$
 \frac{u'}{u} = - x \frac{1 + C + x^2}{1 - C - x^2}\,\cdot
 $$
 Writing this equation in the form
$$
 \frac{d \ln u}{dx} = x + \frac{2 x}{x^2 + C - 1}
 $$
 and integrating we obtain the following general solution to Eq. (\ref{Ricc.Lin_sec1.eq1}):
 \begin{equation}
 \label{Ricc.Lin_sec1.sol1}
 u = K (x^2 + C - 1) {\rm e}^{x^2/2}, \quad C, K = {\rm const.}
 \end{equation}
 \end{exa}

 \begin{exa}
 \label{Ricc.Lin_sec1.ex2}
 If we apply Eqs. (\ref{Ricc.Lin_eq3})-(\ref{Ricc.Lin_eq4}) to Eq. (\ref{gai.eq18})
 from Example \ref{gai.sub2.ex2},
 $$
 y' = x^2 + (x + x^2)\, y +
 \frac{1}{4} \big(2 x + x^2\big)\, y^2,
 $$
 we obtain the following second-order
 linear equation:
 \begin{equation}
 \label{Ricc.Lin_sec1.eq2}
 u'' - (1 + x)\bigg(x + \frac{2}{2x + x^2}\bigg) u' + \frac{1}{4} x^2\big(2 x + x^2\big) u = 0.
 \end{equation}
 Let us integrate this equation. Writing Eq. (\ref{Ricc.Lin_eq2}) in the form
 $$
 \frac{u'}{u} = - R(x) y,
 $$
 substituting here
  $$
  R(x) = \frac{1}{4} \big(2 x + x^2\big)
  $$
   and the expression (\ref{gai.eq19}) for $y,$ we
 get
 $$
 \frac{u'}{u} = - \frac{1}{2} \big(2 x + x^2\big)\frac{1 + \left(C + \int (x + x^2) {\rm e}^{- x^2/2}\,dx \right) {\rm
 e}^{x^2/2}}{1 - \left(C + \int (x + x^2) {\rm e}^{- x^2/2}\,dx \right) {\rm
 e}^{x^2/2}}\,\cdot
 $$
 The integration yields \ $\ln |u| = \ln |K|  + \phi(x),$ where
 $$
 \phi(x) = - \frac{1}{2} \int \big(2 x + x^2\big)\frac{1 +
\left(C + \int (x + x^2) {\rm e}^{- x^2/2}\,dx \right) {\rm
 e}^{x^2/2}}{1 - \left(C + \int (x + x^2) {\rm e}^{- x^2/2}\,dx \right) {\rm
 e}^{x^2/2}} \,dx.
 $$
Hence, the general solution to Eq. (\ref{Ricc.Lin_sec1.eq2}) has the
form
 \begin{equation}
 \label{Ricc.Lin_sec1.sol2}
 u = K\,{\rm e}^{\phi(x)},
 \end{equation}
 where $\phi(x)$ is given above and $K$ is an arbitrary constant.
 \end{exa}

Applying Eqs. (\ref{Ricc.Lin_eq3})-(\ref{Ricc.Lin_eq4}) to Eq.
(\ref{gai.ric3}) and using the solution procedure for Eq.
(\ref{gai.ric3}) described in Section \ref{gai.sub3}, we obtain the
following general result.

 \begin{thm}
 The general solution of the second-order linear equation
 \begin{equation}
 \label{Ricc.Lin_sec1.eq3}
 u'' - \bigg[Q(x)x + \frac{Q'(x) - k P'(x)}{Q(x) - k P(x)}\bigg] u'
 + k \big[P(x) Q(x) - k P^2(x)\big] u = 0
 \end{equation}
 with an arbitrary constant $k$ and two arbitrary functions $P(x)$ and $Q(x)$ can be obtained
 by quadratures.
 \end{thm}

 \section{From second-order to Riccati equations}
 \label{Ricc.Lin_sec2}
 \setcounter{equation}{0}

% \subsection{General approach}

It is manifest from Eqs. (\ref{Ricc.Lin_eq4}) that \textit{two}
 coefficients $f(x)$ and $g(x)$ of a given second-order equation (\ref{Ricc.Lin_eq3}) do not uniquely
 determine \textit{three} coefficients $P(x), Q(x), R(x)$ of the corresponding Riccati equation
 (\ref{Ricc.Lin_eq1}). Namely, if we know solutions of an equation
 (\ref{Ricc.Lin_eq3}), we can solve  by using the formula (\ref{Ricc.Lin_eq2})
 an infinite set of the Riccati equations
 \begin{equation}
 \label{Ricc.Lin_sec2.eq1}
 y' = R(x) y^2 - \bigg[f(x) + \frac{R'(x)}{R(x)}\bigg]\, y + \frac{g(x)}{R(x)}
 \end{equation}
 with an arbitrary function $R(x) \not= 0.$

 \begin{exa}
 \label{Ricc.Lin_sec2.ex1}
 Consider the following equation with constant coefficients:
 \begin{equation}
 \label{Ricc.Lin_sec2.eq2}
  u'' + u = 0.
 \end{equation}
 Here $f = 0, \ g = 1.$ Hence, the corresponding Riccati equation
 (\ref{Ricc.Lin_sec2.eq1}) has the form
 \begin{equation}
 \label{Ricc.Lin_sec2.eq2C}
  y' = R(x)\,y^2 - \frac{R'(x)}{R(x)}\,y + \frac{1}{R(x)}\,\cdot
 \end{equation}
 Substituting the general solution
 $$
 u = C_1 \cos x + C_2 \sin x
 $$
 of Eq. (\ref{Ricc.Lin_sec2.eq2}) in (\ref{Ricc.Lin_eq2}) we
 obtain the following solution to Eq. (\ref{Ricc.Lin_sec2.eq2C}):
 $$
 y = \frac{1}{R(x)}\, \frac{C_1 \sin x - C_2 \cos x}{C_1 \cos x + C_2 \sin
 x}\,\cdot
 $$
 If $C_2 \not= 0$ we denote $K = C_2/C_1$ and writhe the solution in the
 form
 $$
 y = \frac{1}{R(x)}\,\frac{\sin x - K \cos x}{\cos x + K \sin
 x}
 $$
 or, upon dividing the numerator and denominator by $\cos x,$
 \begin{equation}
 \label{Ricc.Lin_sec2.eq2B}
 y = \frac{1}{R(x)}\, \frac{{\rm tg} x - K}{1 + K {\rm tg} x}\,,
 \quad K = {\rm const.}
 \end{equation}
 If $C_2 = 0$ the solution becomes
 $$
 y = - \frac{{\rm ctg} x}{R(x)}
 $$
 which can be obtained from (\ref{Ricc.Lin_sec2.eq2B}) by letting
 $K \rightarrow \infty.$ Thus, the general solution to the Riccati
 equation (\ref{Ricc.Lin_sec2.eq2A}) is given by (\ref{Ricc.Lin_sec2.eq2B})
 where $- \infty \leq K \leq + \infty.$

 In particular, taking in (\ref{Ricc.Lin_sec2.eq2A}), (\ref{Ricc.Lin_sec2.eq2B})
 $R(x) = {\rm e}^x$ we conclude that
 the  equation
 \begin{equation}
 \label{Ricc.Lin_sec2.eq2A}
  y' = {\rm e}^x\,y^2 - y + {\rm e}^{-x}
 \end{equation}
 has the general solution
 \begin{equation}
 \label{Ricc.Lin_sec2.eq2A1}
 y = \frac{{\rm tg} x - K}{1 + K {\rm tg} x}\,{\rm e}^{-x}, \quad
 - \infty \leq K \leq + \infty.
 \end{equation}
 \end{exa}

% \subsection{Particular constructions}

 Sometimes it is convenient to use a
 restricted  correspondence between second-order and
 Riccati equations by writing
 Eqs. (\ref{Ricc.Lin_eq3})-(\ref{Ricc.Lin_eq4}) corresponding to the
 Riccati equation (\ref{Ricc.Lin_eq1}) in
 the following form:
 \begin{equation}
 \label{Ricc.Lin_sec2.eq3}
 R(x)\, u'' - \big[R'(x) + Q(x) R(x)\big] u' + P(x) R^2 (x)\, u =
 0.
 \end{equation}

 \begin{exa}
 \label{Ricc.Lin_sec2.ex2}
 Consider Euler's equation written in the form (\ref{NC.eq12}):
 \begin{equation}
 \label{Ricc.Lin_sec2.eq4}
  x^2\, u'' + (A + 1)\, x\, u' + B\, u = 0, \quad A, B = {\rm const.}
 \end{equation}
 Comparing the equations (\ref{Ricc.Lin_sec2.eq3}) and (\ref{Ricc.Lin_sec2.eq4})
 we take $R(x) = x^2$ and obtain
 $$
 Q(x) = - \frac{A + 3}{x}\, \quad P(x) = \frac{B}{x^4}\,\cdot
 $$
 Thus, we have arrived at the following Riccati equation:
 \begin{equation}
 \label{Ricc.Lin_sec2.eq5}
 y' = x^2 y^2 - \frac{A + 3}{x}\, y + \frac{B}{x^4}\,\cdot
 \end{equation}
 We know that the solutions of Eq. (\ref{Ricc.Lin_sec2.eq4}) have
 the form $u = x^\lambda.$ Substitution in (\ref{Ricc.Lin_eq2})
 yields the following form of solutions to Eq. (\ref{Ricc.Lin_sec2.eq5}):
 \begin{equation}
 \label{Ricc.Lin_sec2.eq6}
 y = - \frac{\lambda}{x^3}\,\cdot
 \end{equation}
 Substituting (\ref{Ricc.Lin_sec2.eq6}) in Eq.
 (\ref{Ricc.Lin_sec2.eq5}) we  obtain again the characteristic
 equation (\ref{NC.eq6}):
  \begin{equation}
 \label{Ricc.Lin_sec2.eq7}
 \lambda^2 + A \lambda + B = 0.
 \end{equation}
However, Eqs. (\ref{Ricc.Lin_sec2.eq6}), (\ref{Ricc.Lin_sec2.eq7})
provide only two particular solutions to Eq.
(\ref{Ricc.Lin_sec2.eq5}). In order to find the general solution of
the \textit{nonlinear} equation (\ref{Ricc.Lin_sec2.eq5}), we have
to construct the general solution $u(x)$ of the \textit{linear}
equation (\ref{Ricc.Lin_sec2.eq4}) using the \textit{superposition
principle} and substitute $u = u(x)$ in (\ref{Ricc.Lin_eq2}).
 \end{exa}

 Using the results of Section \ref{NC.gen} on integrability of Eq.
 (\ref{NC.gen6}),
 $$
  \phi^2\,y'' + \big(A + \phi' - 2 \sigma\big) \phi\, y' +
  \big(B - A \sigma + \sigma^2 - \phi \sigma'\big) y = 0,
  $$
 we can formulate the following general result.

 \begin{thm}
 \label{Ricc.Lin.thm}
 The Riccati equation
 \begin{equation}
 \label{Ricc.Lin_sec2.eq8}
 y' = R(x)\, y^2 - \bigg[\frac{A + \phi\,' - 2 \sigma}{\phi} +
 \frac{R'}{R}\bigg]\, y + \frac{B - A \sigma + \sigma^2 - \phi \sigma\,'}{R \phi^2}
 \end{equation}
 with three arbitrary functions $R(x), ~\phi(x), ~\sigma(x)$ and
 two arbitrary constants $A, ~B$ is integrable by quadratures.
 \end{thm}

 \section{Application to Ermakov's equation}
 \label{Ricc.Erm}
 \setcounter{equation}{0}

 The above results on integration of linear equations
   \begin{equation}
  \label{Ricc.Erm.eq1}
  u'' + a(x) u' + b(x) u = 0
  \end{equation}
 can be combined with Ermakov's method for
 solving nonlinear equations of the\linebreak following form (see Editor's preface to
 Ermakov's paper in this volume):
  \begin{equation}
  \label{Ricc.Erm.eq2}
  u'' + a(x) u' + b(x) u = \frac{\alpha}{u^3}\,{\rm e}^{- 2 \int a(x)
  dx}, \quad \alpha = {\rm const.}
  \end{equation}

 \begin{exa}
 Using the solution (\ref{Ricc.Lin_sec1.sol1}) of Eq. (\ref{Ricc.Lin_sec1.eq1}) and applying Ermakov's method
 one can
 solve the following nonlinear equation:
  \begin{equation}
 \label{Ricc.Erm.eq3}
 u'' - \bigg(2x + \frac{1}{x}\bigg) u' + x^2 u = \alpha x^2 {\rm e}^{2
 x^2}\,u^{-3}.
 \end{equation}
 \end{exa}

 \begin{exa}
 The nonlinear equation  (\ref{Ricc.Erm.eq2}) associated with the integrable equation
 (\ref{NC.gen6}) has the form
  \begin{equation}
 \label{Ricc.Erm.eq4}
 \phi^2\,u'' + \big(A + \phi' - 2 \sigma\big) \phi\, u' +
  \big(B - A \sigma + \sigma^2 - \phi \sigma'\big) u =
  \frac{\alpha}{u^3}\,{\rm e}^{2 \int (2 \sigma - A)/\phi\, dx},
 \end{equation}
 where $\phi$ and $\sigma$ are arbitrary functions of $x,$ and
 $A$ is an arbitrary constant. Eq. (\ref{Ricc.Erm.eq4}) is
 integrable by quadratures.
  \end{exa}

% \begin{thebibliography}{99}


\begin{thebibliography}{10}

 %\begin{thebibliography}{1}

\bibitem{step58}
V.~V. Stepanov, {\em Course of differential equations}.
\newblock Moscow: {S}tate {P}ublisher of {P}hys.-{M}ath. {L}iterature, 1958.
\newblock Russian, 7th ed.

\bibitem{ibr06P}
N.~H. Ibragimov, {\em A practical course in differential equations
and
  mathemtical modelling}.
\newblock Karlskrona: ALGA Publications, 3rd~ed., 2006.

\bibitem{ibr99}
N.~H. Ibragimov, {\em Elementary Lie group analysis and ordinary
differential
  equations}.
\newblock Chichester: John Wiley \& Sons, 1999.

\bibitem{ibr89}
N.~H. Ibragimov, {\em Primer of group analysis}.
\newblock Moscow: Znanie, No. 8, 1989.
\newblock (Russian). Revised edition in English: {\it Introduction to modern
  group analysis}, Tau, Ufa, 2000. Available also in Swedish: {\it Modern
  grouppanalys: En inledning till Lies l\"{o}sningsmetoder av ickelinj\"{a}ra
  differentialekvationer,} Studentlitteratur, Lund, 2002.

\end{thebibliography}
 \end{document}